\documentclass[10pt]{article}
\usepackage{amsfonts,color}
\usepackage{amsmath,amssymb}
\usepackage{epsfig}
\usepackage{lscape}
\usepackage{amssymb}
\usepackage{graphicx}
\usepackage{footnote}
\usepackage[utf8]{inputenc}
\usepackage{mathtools}
\usepackage{amsthm,amssymb,wasysym}
\usepackage{amsfonts}
\usepackage[english]{babel}
\usepackage{tikz}
\usepackage{bbm}
\usepackage{colonequals}
\usepackage{setspace}
\usepackage{changes}

\allowdisplaybreaks[1]

\makeindex

\def\12{{\frac{1}{2}}}                             
\newcommand{\f}[1]{{\bf{#1}}}                      
\newcommand{\fg}[1]{{\mbox{\boldmath $#1$}}}       
\newtheorem{theorem}{Theorem}[section]

\newtheorem{remark}[theorem]{Remark}

\newcommand{\R}{\mathbb{R}}

\newcommand{\pt}{\partial}

\newcommand{\norm}[1]{\|#1\|}                      
\newcommand{\m}{\mathfrak{m}}                      

\DeclareMathOperator{\sym}{sym}
\DeclareMathOperator{\tr}{tr}

\DeclareMathOperator{\axl}{axl}

\DeclareMathOperator{\anti}{anti}
\DeclareMathOperator{\dev}{dev}
\DeclareMathOperator{\sL}{\mathfrak{sl}}
\DeclareMathOperator{\so}{\mathfrak{so}}

\DeclareMathOperator{\Lin}{Lin}

\DeclareMathOperator{\Curl}{Curl\,}
\DeclareMathOperator{\curl}{curl\,}
\newcommand{\Sym}{ {\rm{Sym}} }

\newcommand{\id}{{\boldsymbol{\mathbbm{1}}}}

\def\div{\textrm{div}}
\def\Div{\textrm{Div\,}}
\def\DIV{\textrm{DIV\,}}
\newcommand{\GRAD}[1]{{\rm GRAD}[{#1}]}
\newcommand{\Grad}[1]{{\rm Grad}[{#1}]}
\newcommand{\grad}[1]{{\rm grad}[{#1}]}

\newcommand{\ks}{\widetilde{\f k}}        
\newcommand{\sigmas}{\widetilde{\sigma}}        
\newcommand{\xbar}{\bar{\bf x}}            

\def\skew{\text{skew} \, }

\newcommand{\di}{\,{\rm d}}

\newcommand{\Scal}[1]{\langle#1\rangle}            

\DeclareMathOperator{\Euklid}{\mathbb{R}^3}        
\DeclareMathOperator{\Skalar}{\mathbb{R}}          
\newcommand{\B}{{\cal B}}                          

\makeatletter
\let\@fnsymbol\@arabic

\newcommand{\TD}{\mathbb{D} \widetilde{{\mbox{\boldmath $\sigma$}}}({\bf x}_0)}
\newcommand{\TDnp}{\mathbb{D}^{\rm np} \widetilde{{\mbox{\boldmath $\sigma$}}}({\bf x}_0)}
\newcommand{\TDp}{\mathbb{D}^{\rm p} \widetilde{{\mbox{\boldmath $\sigma$}}}({\bf x}_0)}
\newcommand{\TDpone}{\mathbb{D}^{\rm p1} \widetilde{{\mbox{\boldmath $\sigma$}}}({\bf x}_0)}
\newcommand{\TDptwo}{\mathbb{D}^{\rm p2} \widetilde{{\mbox{\boldmath $\sigma$}}}({\bf x}_0)}
\newcommand{\Deltax}{\Delta {\bf x}}

\newcommand{\TDb}{\mathbb{D}_2^{\rm b} \widetilde{{\mbox{\boldmath $\sigma$}}}({\bf x}_0)}
\newcommand{\TDbone}{\mathbb{D}_{2}^{\rm b1} \widetilde{{\mbox{\boldmath $\sigma$}}}({\bf x}_0)}
\newcommand{\TDbtwo}{\mathbb{D}_{2}^{\rm b2} \widetilde{{\mbox{\boldmath $\sigma$}}}({\bf x}_0)}

\newcommand{\Deltaxb}{\Delta {\bf x}^2_{\rm b}}
\newcommand{\TDq}{\mathbb{D}_2^{\rm q} \widetilde{{\mbox{\boldmath $\sigma$}}}({\bf x}_0)}
\newcommand{\Deltaxq}{\Delta {\bf x}^2_{\rm q}}

\begin{document}
\title{The modified indeterminate couple stress model: Why Yang et al.'s arguments motivating a symmetric couple stress tensor contain a gap and why the couple stress tensor may be chosen symmetric nevertheless
}
\author{\normalsize{Ingo M\"unch\thanks{Corresponding author: Ingo M\"unch, Institute for Structural Analysis, Karlsruhe Institute of
Technology, Kaiserstr. 12, 76131 Karlsruhe, Germany, email: ingo.muench@kit.edu}
\quad and \quad
Patrizio Neff\thanks{Patrizio Neff,  \ \ Head of Lehrstuhl f\"{u}r Nichtlineare
Analysis und Modellierung, Fakult\"{a}t f\"{u}r Mathematik, Universit\"{a}t Duisburg-Essen,
Thea-Leymann Str. 9, 45127 Essen, Germany, email: patrizio.neff@uni-due.de}
\quad and \quad
Angela Madeo\footnote{Angela Madeo, \ \  Laboratoire de G\'{e}nie Civil et Ing\'{e}nierie
Environnementale, Universit\'{e} de Lyon-INSA, B\^{a}timent Coulomb, 69621 Villeurbanne Cedex, France;
and International Center M\&MOCS ``Mathematics and Mechanics of Complex Systems", Palazzo Caetani,
Cisterna di Latina, Italy,
 email: angela.madeo@insa-lyon.fr}
\quad and \quad
Ionel-Dumitrel Ghiba\thanks{Ionel-Dumitrel Ghiba, \ \ \ \ Lehrstuhl f\"{u}r Nichtlineare Analysis und
Modellierung, Fakult\"{a}t f\"{u}r Mathematik, Universit\"{a}t Duisburg-Essen, Thea-Leymann Str. 9,
45127 Essen, Germany;  Alexandru Ioan Cuza University of Ia\c si, Department of Mathematics,  Blvd.
Carol I, no. 11, 700506 Ia\c si, Romania; and  Octav Mayer Institute of Mathematics of the Romanian
Academy, Ia\c si Branch,  700505 Ia\c si, email: dumitrel.ghiba@uni-due.de, dumitrel.ghiba@uaic.ro} }
} \maketitle


\begin{abstract}
We show that the reasoning in favor of a symmetric couple stress tensor in Yang et al.'s introduction of the modified couple stress theory contains a gap, but we present a reasonable physical hypothesis, implying that the couple stress tensor is traceless and may be symmetric anyway. To this aim, the origin of couple stress is discussed on the basis of certain properties of the total stress itself. In contrast to classical continuum mechanics, the balance of linear momentum and the balance of angular momentum are formulated at an infinitesimal cube considering the total stress as linear and quadratic approximation of a spatial Taylor series expansion.
\\
\vspace*{0.25cm}

\noindent AMS Math 74A10 (Stress), 74A35 (Polar materials), 74B05 (Classical linear elasticity)
\end{abstract}

\tableofcontents

\section{Introduction}

Toupin~\cite{Toupin62,Toupin64} and Mindlin et al.
\cite{Mindlin62,Mindlin63,Mindlin64,Mindlin65,Mindlin68} have established strain gradient theories to evaluate size effects by allowing the total stress tensor to become asymmetric. Additionally, the work conjugated quantity to the strain gradient, namely the couple stress tensor, was also accepted to be asymmetric. Such models usually reproduce the size effect in the sense that "smaller is stronger", which is a central point of strain gradient modeling \cite{Fleck93,Fleck94,Fleck96,Fleck97}.

The special strain gradient theory in the work of Yang et al.~\cite{Yang02} is a subclass of the former and uses a symmetric second order couple stress tensor ${\f m}$ for the so called modified couple stress model, whose decisive advantage is to reduce the number of additional constitutive coefficients to just one characteristic length scale. In effect, they try to motivate that the couple stress tensor itself should be symmetric. Many subsequent papers have used this approach. From our point of view, an artificial equilibrium condition is used to imply the symmetry of the couple stress tensor. Therefore, their argument is not consistent, as has also been previously noted by Lazopoulos~\cite{Lazo09} and Hadjesfandiari and Dargush~\cite{hadjesfandiari2014evo}. However, this does not mean that symmetry of the couple stress tensor in itself violates any physical law, as repeatedly claimed by Hadjesfandiari and Dargush~\cite{hadjesfandiari2014evo}. Indeed, several different motivations for a symmetric couple stress tensor have already been given. For example, a sequence of smaller and smaller samples should have bounded stiffness in bending and torsion since the physics dictates bounded energy. However, in the framework of the isotropic, linear Cosserat model (and the indeterminate couple stress model is a limit of that model) such boundedness necessitates to take a symmetric moment stress tensor  ${\f m}$ \cite{Neff_Jeong_Conformal_ZAMM08,Jeong_Neff_ZAMM08}.\\

In this work we will therefore review the indeterminate couple stress model in order to
appreciate some of the modeling issues which are involved in the discussion above. First, we are presenting the couple stress framework based on a variational derivation. This derivation starts from assuming a certain energy involving second (rotational) gradients and immediately uses the small strain kinematics and the isotropy assumption. In this way, a format of the balance equations is derived as Euler-Lagrange equations, together with the assumed constitutive relations and various boundary conditions. Here, the couple stress tensor is already seen to be trace-free: $\tr({\f m}) = 0$. In addition we touch further upon conformal invariance requirements which are naturally connected with symmetric couple stresses.\\

This variational development will be compared with another derivation of the couple stress model, which does not make use of any constitutive assumptions like small strains, linearity or isotropy. We will only invoke balance principles applied to infinitesimal cubes and the existence of a total (in general non-symmetric) stress field $\widetilde{\fg \sigma}$. The couple stress tensor ${\f m}$ will be identified with its assignment to be an exclusive stress resultant for the balance of angular momentum. On the other hand, the couple stress tensor will not appear in the balance of linear momentum.

However, contrarily to what is classically done, we allow higher order variations of the total stress field over the cube with the result of being able to clearly identify the couple stress tensor $\f m$. The final outcome is a set of two coupled balance equations having exactly the same format as the indeterminate couple stress model derived with the variational approach. Since no constitutive assumptions are yet involved, no condition for the trace of the couple stress tensor is included. Interestingly, if we assume from the outset the symmetry of the total stress tensor $\widetilde{\fg \sigma}$ we arrive consistently at $\tr({\f m}) = 0$ in this framework. However, the total stress tensor in the indeterminate couple stress model is not symmetric in the standard approach. This can be modified by adding a self-equilibrated stress-field (thus, the balance equations are the same, but different boundary conditions are implied) and the full correspondence can be established to the model, see our paper \cite{GhiNeMaMue15}. Additionally, relating couple stress effects to isochoric deformation modes only, we find again an argument for the symmetry of the couple stress tensor $\f m$. Both results can also be found in continuum theories with conformally invariant curvature measure, which reduces the number of constitutive parameters to a minimum
\cite{Neff_Jeong_Conformal_ZAMM08, Neff_Jeong_bounded_stiffness09}. Thus, our procedure provides constitutive statements from an equilibrium method and reasonable assumptions of the classical continuum theory, recovering previous results obtained via variational procedures now  using a suitable balance approach.\\

Finally, we critically discuss in detail the underlying reasoning of the motivation of a symmetric couple stress tensor in Yang et al.~\cite{Yang02}  and come to the conclusion that it is not tenable: the authors do not present a convincing argument for symmetric couple stress tensors. This, however, does not imply that assuming symmetry of the couple stress tensor violates any fundamental physical law, as erroneously claimed in Hadjesfandiari and Dargush \cite{hadjesfandiari2014evo}. We show this by some simple examples.\\

The paper is structured as follows. The classical indeterminate couples stress model
including some of its variants in the isotropic and hyperelastic setting and some remarks on conformal invariance of the curvature energy are recalled. Then, the balance equations for couple stress models are treated in general. For ease of understanding we deliberately repeat a few formulas in Section \ref{KapIndeterminate} and \ref{KapRelatedModels}. Next, the origin and properties of couple stress are systematically investigated on the basis of the total stress function, the balance of linear momentum, and the balance of angular momentum. Then we try to explain the approach of Yang et al.~\cite{Yang02} bona fide and indicate where their argument fails. Finally, we give an analytical example to verify results of this paper.
\subsection{Notational agreements}\label{KapNotations}
With $\R^{3\times 3}$ we denote the set of real $3\times 3$ second order tensors, written with bold capital letters. Vectors in $\R^{3}$ are denoted by small bold letters. The
components of tensors and vectors are given according to orthogonal unit vectors $\f e_1, \, \f e_2, \, \f e_3$. We use Lagrangian coordinates $\f x$ to describe physical fields.
Throughout this paper (when not specified otherwise) Latin subscripts specify the
direction of components and take the values $1,2,3$. For repeating subscripts Einstein's summation convention applies.

For vectors $\f a, \f b\in\R^3$ we let $\langle {\f a},{\f b}\rangle_{\R^3}$ denote the canonical scalar product on $\R^3$ with associated vector norm $\|{\f a}\|^2_{\R^3}=\langle {\f a},{\f a}\rangle_{\R^3}$. The standard Euclidean scalar product on $\R^{3\times 3}$ is given by $\langle{\f X},{\f Y}\rangle_{\R^{3\times3}}=\tr({\f X \f Y^T})$, and thus the Frobenius tensor norm is $\|{\f X}\|^2=\langle{\f X},{\f X}\rangle_{\R^{3\times3}}$. The identity tensor on $\R^{3\times3}$ will be denoted by $\id$, so that $\tr({\f X})=\langle{\f X},{\id}\rangle$. We adopt the following abbreviations:  $\so(3)\colonequals\{\f X\in\mathbb{R}^{3\times3}\;|\f X^T=- \f X\}$ is the vector space of  skew-symmetric tensors and $\sL(3)\colonequals\{\f X\in\mathbb{R}^{3\times3}\;| \tr({\f X})=0\}$ is the vector space of traceless tensors. For all $\f X\in\mathbb{R}^{3\times3}$ we set $\sym \f X=\frac{1}{2}(\f X^T+\f X)\in\Sym(3)$,
$\skew \f X=\frac{1}{2}(\f X-\f X^T)\in \so(3)$ and the deviatoric part $\dev \f X=\f
X-\frac{1}{3}\;\tr(\f X)\cdot \,\id\in \sL(3)$  and we have the decomposition
of $\Skalar^{3 \times 3}$
\begin{align}
\Skalar^{3 \times 3}=\{\sL(3)\cap \Sym(3)\}\oplus\so(3) \oplus\mathbb{R}\!\cdot\! \id \, ,\qquad \f X=\dev \sym \f X+ \skew \f X+\frac{1}{3}\tr(\f X) \, \id \, ,
\end{align}
simply allowing to split every second order tensor $\f X\in\mathbb{R}^{3\times3}$ uniquely into its trace free symmetric part, skew-symmetric part and spherical part, respectively. Typical conventions for differential operations are implied such as comma followed by a subscript to denote the partial derivative with respect to  the corresponding cartesian coordinate. The outer product of vectors
\begin{align}\label{crossProd}
\f a \times \f b = a_i \, b_j \, \epsilon_{ijk} \, \f e_k
= \left( \begin{array}{c}
a_2 \, b_3 - a_3 \, b_2 \\
a_3 \, b_1 - a_1 \, b_3 \\
a_1 \, b_2 - a_2 \, b_1 \\
\end{array}\right) \, ,
\end{align}
is given by the components of the alternating Levi-Civit\`{a} tensor $\fg \epsilon$ with $\epsilon_{ijk}=+1$ for $ijk=\{123,231,312\}$, $\epsilon_{ijk}=-1$ for $ijk=\{132,213,321\}$, and $\epsilon_{ijk}=0$ else. Let us define the operator $\axl:\so(3)\rightarrow\mathbb{R}^3$ such that
\begin{align}\label{DefaxlMitTimes}
\f A \cdot \f b = \axl [\f A] \times \f b \, .
\end{align}
In symbolic and index notation the axial vector $\f a$ reads
\begin{align}\label{Defaxl}
\f a = \axl [\f A] \colonequals -\frac{1}{2} \, \f A : \fg \epsilon \, \Leftrightarrow \, a_k \colonequals -\frac{1}{2} \, A_{ij} \, \epsilon_{ijk} \, ,
\end{align}
where the colon in $\f A : \fg \epsilon $ denotes double contraction. The inverse function of $\axl$ can be found by using eq.\eqref{Defaxl} in index notation multiplied by $\epsilon_{abk}$ yielding
\begin{align}\label{Defanti}
& \epsilon_{abk} \, a_k = -\frac{1}{2} \, A_{ij} \, \epsilon_{ijk} \,\epsilon_{abk} \, \Leftrightarrow \,
\epsilon_{abk} \, a_k = -\frac{1}{2} \, A_{ij} \, (\delta_{ia} \, \delta_{jb} - \delta_{ib} \, \delta_{ja})  \notag \\
& \Leftrightarrow \, \epsilon_{abk} \, a_k = -\frac{1}{2} \, (A_{ab} - A_{ba}) \,  \Leftrightarrow \, \epsilon_{abk} \, a_k = - A_{ab} \, \Leftrightarrow \, \epsilon_{abk} \, a_k = -\frac{1}{2} \, (A_{ab} - A_{ba}) \notag \\
& \Leftrightarrow \, \f A = \anti [\f a] \colonequals \, - \fg \epsilon \cdot \f a \, \Leftrightarrow \, A_{ab} \colonequals -\epsilon_{abk} \, a_k \, ,
\end{align}
since $\f A = - \f A^T$. Thus, the components of the axial vector $\f a$ define the components of the skew-symmetric tensor $\f A$ according to
\begin{align}
\f a=\left(\begin{array}{c}
a_1\\
a_2\\
a_3\end{array}\right) \, , \,
\f A=\left(\begin{array}{ccc}
0 &-a_3&a_2\\
a_3&0& -a_1\\
-a_2& a_1&0
\end{array}\right) \, .
\end{align}
Gradients of scalar fields $\phi \in \Skalar$, vector fields $\f b \in \Euklid$, and second order tensor fields $\f X\in\mathbb{R}^{3\times3}$ are defined by
\begin{align}\label{GradSkalDef}
\grad{\phi}  \colonequals \frac{\partial \, \phi}{\partial \, \f x} = \frac{\partial \, \phi}{\partial \, x_i} \, \f e_i = \phi_{,i} \, \f e_i \quad \in \Euklid \,,
\end{align}
\begin{align}\label{GradVektDef}
\Grad{\f b}  \colonequals \frac{\partial \, \f b}{\partial \, \f x} = \frac{\partial \, b_i}{\partial \, x_j} \, \f e_i \otimes \f e_j = b_{i,j} \, \f e_i \otimes \f e_j \quad \in \Skalar^{3 \times 3} \,,
\end{align}
\begin{align}\label{GradTensor2Def}
\GRAD{\f X}  \colonequals \frac{\partial \, \f X}{\partial \, \f x} = \frac{\partial \, X_{ij}}{\partial \, x_k} \, \f e_i \otimes \f e_j \otimes \f e_k = X_{ij,k} \, \f e_i \otimes \f e_j \otimes \f e_k \quad \in \Skalar^{3 \times 3 \times 3} \,.
\end{align}
The divergence of vector fields, second order tensor fields, and third order tensor fields reads
\begin{align}\label{VektDiv}
\div \f b  \colonequals \tr ( \Grad{\f b}) = \langle \Grad{\f b} , \id \rangle = b_{i,j} \, \delta_{ij} =  b_{i,i} \quad \in \Skalar \, ,
\end{align}
\begin{align}\label{TensorDiv}
\Div \f X  \colonequals (\GRAD{\f X})_{ijk} \, \delta_{jk} \, \f e_i = X_{ij,j} \, \f e_i \quad \in \Euklid \, ,
\end{align}
\begin{align}\label{TensorDivThirdOrder}
\DIV \f m  \colonequals
m_{ijk,k} \, \f e_i \otimes \f e_j \quad \in \Skalar^{3 \times 3} \, .
\end{align}
Finally, we define the curl of vector fields and of second order tensor fields according to
\begin{align}\label{CurlDefVektor}
\curl \f v \colonequals  -  v_{a,b} \, \epsilon_{abi} \, \f e_i \quad \in \Euklid \, ,
\end{align}
\begin{align}\label{CurlDefTensor}
\Curl {\f X} \colonequals  - X_{ia,b} \, \epsilon_{abj} \, \f e_i \otimes \f e_j \quad \in \Skalar^{3 \times 3} \, .
\end{align}
Using eq.\eqref{Defaxl}  the curl of vectors can be written in terms of gradients by
\begin{align}\label{axlgradvcurl}
\curl \f v = 2 \, \axl ( \skew \Grad{\f v}) & = - \, (\Grad{\f v})_{ab} \, \epsilon_{abk} \, \f e_k \,.
\end{align}
We consider an elastic body $\B$ occupying a bounded open set of the three-dimensional Euclidian space $\R^3$ with boundary $\partial \B$ as piecewise smooth surface. In $\B \subset \R^3$ submerged subdomains $\B_c$ exist, such that their surfaces $\partial \B_c$ do not overlap any portion of $\partial \B$. In our notation the Cauchy theorem reads
\begin{align}\label{CauchyFundaTheoremI}
\f t(\f x) = \fg \sigma(\f x) \cdot \f n(\f x)  \, , \quad \f t, \, \f n \in \R^3 , \, \fg \sigma \in \R^{3 \times 3} \, ,
\end{align}
where the traction $\f t$ at $\f x$ is given by the stress tensor $\fg \sigma(\f x)$ and the surface normal $\f n(\f x)$. In correspondence with the notation in Yang et al. \cite{Yang02}, we consider this notation, which is common in American literature, e.g., Truesdell \& Noll \cite{Truesdell65}, Marsden \& Hughes \cite{MarsdenHughes1994}, and Gurtin \cite{Gurtin2003}. However, different invariant representations of tensor calculus exist and our notation is different to the notation of the tensor calculus introduced by Gibbs. The latter notation is sometimes used in European literature, e.g., Altenbach and coworkers \cite{Altenbach2003, Altenbach2012,Altenbach2004}, Eremeyev \cite{Eremeyev2014}, and Lurie \cite{Lurie2012}. In this tensor calculus the Cauchy theorem reads $\f t(\f x) = \f n(\f x) \cdot \fg \sigma(\f x)$. The difference between both notations results from the transposition in the definition of the stress tensor together with the transposed representation of the deformation gradient.
\section{The isotropic linear indeterminate couple stress theory}\label{KapIndeterminate}\setcounter{equation}{0}
In this paper we limit our analysis to frame indifferent models with isotropic material, and only to the second gradient of the displacement:
\begin{align}\label{D2u}
{\rm D}^2_{\f x} \f u=\underbrace{\frac{\partial^2 u_i}{\partial x_j \, \partial
x_k}}_{u_{i,jk}} \, \f e_i \otimes \f e_j \otimes \f e_k =(\varepsilon_{ji,k}+\varepsilon_{ki,j}-\varepsilon_{jk,i}) \,\f e_i \otimes \f e_j \otimes \f e_k \, ,
 \end{align}
where
\begin{align}\label{strain}
\fg \varepsilon= \varepsilon_{ij} \, \f e_i \otimes \f e_j = \12 (u_{i,j} +
u_{j,i})\, \f e_i \otimes \f e_j = \12 (\text{Grad} [\f u ]+ (\text{Grad} [\f u])^T ) = \sym \text{Grad}[ \f u] \, .
 \end{align}
is the symmetric linear strain tensor. Thus, from eq.\eqref{D2u} all second derivatives ${\rm D}^2_{\f x} \f u$ of the displacement field $\f u$ can be obtained from linear combinations of $\Grad{\fg \varepsilon}$. In general, strain gradient models do not introduce additional independent degrees of freedoms\footnote{In contrast to Cosserat models where an independent rotation field is under consideration.} aside the displacement field $\f u$. Thus, the higher derivatives introduce a ``latent-microstructure" (constraint microstructure \cite{Grioli03}). However, this apparent simplicity has to be payed with more complicated and intransparent boundary
conditions, as treated in a series of papers \cite{MaGhiNeMue15, NeGhiMaMue15, NeMueGhiMa15}.\\

The linear indeterminate couple stress model is a particular second gradient elasticity model, in which the higher order interaction via moment stresses is restricted to the gradient of the continuum rotation $\curl \, \f u$, where $\f u : \B \mapsto \R^3$ is the displacement of the body. The linear indeterminate couple stress model is therefore interpreted  to be sensitive to rotations of material points and it is possible to prescribe boundary conditions of rotational type. Superficially, this is the simplest possible generalization of linear elasticity in order to include the gradient of the local continuum rotation as a source of additional strains and stress with an associated energy.\\

\noindent
Further, we assume the isotropic quadratic elastic energy to be given by
\begin{align}\label{Welastic}
W(\fg \varepsilon , \widetilde{\f k})&= \int_\B W_{\rm lin}(\fg \varepsilon) + W_{\rm curv}(\ks) \, \di V \notag\\
&= \int_\B \mu \norm{\fg \varepsilon}^2 + \frac{\lambda}{2} [\tr(\fg \varepsilon)]^2 + \mu\,L_c^2\,(\alpha_1\, \|\dev\sym [\widetilde{\f k}]\|^2 + \alpha_2\, \|\skew[\widetilde{\f k}]\|^2) + \alpha_3 \, [\tr(\widetilde{\f k})]^2  \, \di V \, ,
\end{align}
where $\mu$ and $\lambda$ are the classical Lam\'{e} constants and the curvature energy
is expressed in terms of the second order curvature tensor
\begin{align}\label{kcurlu}
\widetilde{\f k} \colonequals \Grad{ \axl(\skew \Grad{ \f u})}=\12 \, \Grad{ \curl \, \f u} \, ,
\end{align}
with additional dimensionless constitutive parameters $\alpha_1$, $\alpha_2$, $\alpha_3$, and $L_c > 0$ as characteristic length. Taking free variations $\delta \f u\in  C^2(\overline{\Omega})$ of the elastic energy $W(\fg \varepsilon,\widetilde{\f k})$ yields the virtual work principle
\begin{align}\label{gradeq211}
\frac{\rm d}{\rm dt}W(\Grad{ \f u}+t \, \Grad{ \delta \f u})=&\int_\B 2\mu\,\langle \fg \varepsilon, \Grad{ \delta \f u} \rangle+\lambda \tr(\fg \varepsilon)\,\tr( \Grad{ \delta \f
u)}\notag\\
&+ 2 \, \mu \, L_c^2 \, \alpha_1\, \langle \dev \sym (\ks),\dev\sym \Grad{\axl \, \skew \Grad{ \delta \f u}} \rangle \notag \\
& + 2 \, \mu \, L_c^2 \, \alpha_2\, \langle \skew(\ks) ,
\skew(\Grad{ \axl\,\skew \, \Grad{ \delta \f u}}) \rangle \notag \\
& + 2 \, \mu \, L_c^2 \,  \alpha_3 \, \tr(\ks) \, \tr(\Grad{\axl \, \skew \Grad{ \delta \f u}}) + \langle \f f , \delta \f u\rangle \di
V =0 \, .
\end{align}
Using the classical divergence theorem for the curvature term in eq.\eqref{gradeq211} it
follows after some simple algebra that
\begin{align}\label{germaneq311}
\int_\B \langle \Div ({\fg \sigma}+\widetilde{\fg \tau})+\f f, \delta \f u \rangle \, \di V-\!\int_{\partial \B}\langle ({\fg \sigma}+\widetilde{\fg \tau}) \cdot \,\f  n, \delta \f u\rangle \, \di A + \! \int_{\partial\B}&\langle  \f m \cdot \f n, \axl\skew \Grad{ \delta \f u} \rangle \di A=0 \, ,
\end{align}
where  $\fg \sigma$ is the symmetric local force-stress tensor from isotropic, linear elasticity
\begin{align}
\fg \sigma=2\, \mu \, \fg \varepsilon+\lambda \, \tr( \fg \varepsilon)\, \id \quad \in
{\rm Sym}(3) \,,
\end{align}
and $\widetilde{\fg \tau}$ represents  the additional non-local force-stress  tensor
\begin{align}\label{nonlocalForcestress}
 \widetilde{\fg \tau}&= - \frac{1}{2}\anti {\rm Div}[ {\f m}] \, \in \so(3) \, ,
\end{align}
which here is automatically skew-symmetric. Using eq.\eqref{Defaxl} we obtain from eq.\eqref{nonlocalForcestress}
\begin{align}\label{Divmplus2axltau}
{\rm Div}[\f m] = - 2 \, \axl [\widetilde{\fg \tau}] \, \Leftrightarrow  {\rm Div}[\f m] + 2 \, \axl [\widetilde{\fg \tau}] = 0 \, .
\end{align}
The second order couple stress\footnote{Also
denominated hyperstress or moment stress.} tensor $\f m$  in eq.\eqref{nonlocalForcestress} reads
\begin{align}\label{DefCoupleStress}
{\f m}&=\mu \, L_c^2 \, [{\alpha_1}\dev \sym (\Grad{ \curl \, \f u}) + {\alpha_2}\,\skew(\Grad{ \curl \, \f u}) + \alpha_3 \, \underbrace{\tr(\Grad{ \curl \, \f u})}_{=0} \, \id] \notag \\
&=2 \, \mu \, L_c^2\,[\alpha_1 \, \dev \sym (\ks) + \alpha_2 \, \skew(\ks)] \,,
\end{align}
which may or may not be symmetric, depending on the material parameters $\alpha_1$, $\alpha_2$. Moreover, $\f m$ in eq.\eqref{DefCoupleStress} is automatically trace free since both the deviator and the skew operator yield trace free tensors.\footnote{One may simplify all formulae using $ \dev \sym (\ks) = \sym (\ks) $.} This is in accordance with our subsequent discussion at an infinitesimal cube in Section \ref{KapCouples}, where symmetric total force stress will yield $\tr (\f m) = 0$ in eq.\eqref{tr_m_from_Sigma}.

Note, however, that the skew-symmetry of the non-local force stress $\widetilde{\fg \tau}$ appears as a constitutive assumption. Thus, if the test function $\delta \f u\in C^{\infty}_0 (\overline{\Omega})$ also satisfies $\axl(\skew \Grad{ \delta \f u})=0$ on $\overline{\Omega}$ (equivalently $\curl \, \delta \f u=0$), then we obtain the balance of momentum
\begin{align}\label{Div2muvarepsilon}
\Div \bigg\{&\underbrace{2\, \mu \, \fg \varepsilon+\lambda \, \tr(\fg \varepsilon)\,
\id}_{\text{\rm local force stress}\  \sigma \, \in \, {\rm Sym}(3)}
\underbrace{- \anti {\rm Div}\{\underbrace{\mu\,L_c^2\,\alpha_1\, \dev\sym
(\ks)+ \mu\,L_c^2\,\alpha_2\,\skew (\ks)}_{\text{\rm hyperstress } \tfrac{1}{2}\f m \, \in \, \Skalar^{3 \times 3}}\}}_{\text{\rm completely skew-symmetric non-local force stress}\ \widetilde{\fg \tau} \, \in \, \so(3)}\bigg\}+ \, \f f=0 \, .
\end{align}
The balance of angular momentum is given by eq.\eqref{nonlocalForcestress}. Both are combined to the compact equilibrium equation
 \begin{align}\label{CoupleStressBalanceMomentum}
\fbox{
\parbox[][2.5cm][c]{14cm}{$
 \widetilde{\fg \sigma} = \fg \sigma +  \widetilde{\fg \tau} \qquad \text {total force stress} \, , \\
 \\
 \Div \widetilde{\fg \sigma} + \f f = 0  \quad \Leftrightarrow \quad
 \Div [\fg \sigma - \frac{1}{2} \anti \Div \f m] + \f f = 0
 \, , \\
 \\
 \Div{\f m} + 2 \axl(\widetilde{\fg \tau}) = 0  \quad \Leftrightarrow \quad
 \Div{\f m} + 2 \axl(\skew \widetilde{\fg \sigma}) = 0 \, .$}}
 \end{align}
A linear Cosserat model with rotation vector $\fg \theta = \axl{\skew{\Grad{\f u}}}$ descents into the indeterminate couple stress model with $\ks =\Grad{\f A(\fg \theta)} = \Grad{\axl{\skew{\Grad{\f u}}}}$.
\subsection{Related models in isotropic second gradient elasticity}\label{KapRelatedModels}
Let us consider the following strain and curvature energy as a minimization problem
\begin{align}
I(\f u)=\int_\B \left[\mu\, \|{\rm sym} \, \Grad{ \f u}\|^2+\frac{\lambda}{2}\, [\tr({\rm sym} \, \Grad{ \f u})]^2+W_{\rm curv}({\rm D}^2_{\f x} \f u)\right] \di V \, \mapsto \,
\text{min. w.r.t.} \, \f u \, ,
\end{align}
admitting unique minimizers under some appropriate boundary condition. Here $\lambda,\mu$ are the Lam\'{e} constitutive coefficients of isotropic linear elasticity, which is
fundamental to small deformation gradient elasticity. If the curvature energy has the
form $W_{\rm curv}({\rm D}^2_{\f x} \f u)=W_{\rm curv}({\rm D}_{\f x} \sym \nabla \f u)$, the model is called {\bf a strain gradient model}. We define the hyperstress tensor of third order as $\m={\rm D}_{{\rm D}^2_{\f x} \f u}W_{\rm curv}({\rm D}^2_{\f x} \f u)$. Since $\Div \m$ is generally not symmetric, the total force stress tensor in a general gradient elasticity theory is not anymore symmetric.

In the following we recall  some curvature energies proposed in different isotropic second gradient elasticity models for the convenience of the reader:
\begin{itemize}
\item {\bf the indeterminate couple stress model} (Grioli-Koiter-Mindlin-Toupin model) \cite{Grioli60,Aero61,Koiter64,Mindlin62,Toupin64,Sokolowski72,Grioli03}
in which the higher derivatives (apparently)  appear only through derivatives of the infinitesimal continuum rotation $\curl \, \f u$.  Hence, the curvature energy  has  the equivalent forms
\begin{align}\label{KMTe}\notag
W_{\rm curv}(\ks)&=\frac{\mu \,L_c^2}{4} \, (\alpha_1 \, \|\sym \Grad{ \curl\, \f u}\|^2+\alpha_2\,\| \skew \Grad{ \curl\, \f u}\|^2\notag\\
&=\mu\,L_c^2 \, (\alpha_1\, \|\sym\underbrace{\Grad{\axl(\skew \Grad{ \f u})}}_{\ks}\|^2 \notag \\
&+\,\alpha_2\,\| \skew \underbrace{ \Grad{\axl(\skew
\Grad{ \f u})}}_{\ks}\|^2 \, , \\
\f m &= 2\,\mu\,L_c^2 \, (\alpha_1\,\sym\ks+\alpha_2\,\skew\ks)  \, .
\end{align}
We remark that the spherical part of the couple stress tensor is zero since $\tr(2\,\ks)=\tr(\nabla \curl \, \f u)={\rm div} (\curl \, \f u)=0$, as seen before. In order to prove the pointwise uniform positive definiteness it is assumed that $\alpha_1>0, \alpha_2>0$. Pointwise uniform positivity is often assumed \cite{Koiter64} when deriving analytical solutions for simple boundary value problems because it allows to invert the couple stress-curvature relation.
\item
  {\bf the modified symmetric couple stress model - the conformal model}.  On the other hand, in the  conformal case \cite{Neff_Jeong_IJSS09,Neff_Paris_Maugin09} one may consider $\alpha_2=0$, which makes  the couple stress tensor ${\f m}$ symmetric and trace free \cite{dahler1963theory}.  This conformal  curvature case has been derived by Neff in \cite{Neff_Jeong_IJSS09}, the curvature energy having the form
\begin{align}
W_{\rm curv}(\ks)&=\mu \, L_c^2\, \alpha_1\, \|\dev \sym \ks\|^2 , \qquad \f m=2\,\mu\,L_c^2\,\alpha_1\,\dev \sym\,\ks\,.
\end{align}
Indeed, there are two major reasons uncovered in \cite{Neff_Jeong_IJSS09} for using the modified couple stress model. First, in order to avoid non-physical singular stiffening behaviour for smaller and smaller samples in bending \cite{Neff_Jeong_bounded_stiffness09} one has to take $\alpha_2=0$. Second, a homogenization procedure invoking a natural ``micro-randomness" assumption (a strong statement of microstructural isotropy) implies conformal invariance, which is again $\alpha_2=0$. Such a model is still well-posed \cite{Neff_JeongMMS08} leading to existence and uniqueness results with only one additional material length scale parameter, although it is {\bf not} pointwise uniformly positive definite.
Since $\axl \skew \widetilde{\f k} = 1/4 \, \curl \curl \f u$ the symmetric couple stress $\f m \in \Sym(3)$ does not work on $\curl \curl \f u$.
\item {\bf the skew-symmetric couple stress model}.
  { Hadjesfandiari and Dargush} strongly advocate
  \cite{hadjesfandiari2011couple,hadjesfandiari2013fundamental,hadjesfandiari2013skew,hadjesfandiari2014evo} the opposite extreme case, $\alpha_1=0$ and $\alpha_2>0$, i.e.~the curvature  energy
\begin{align}
W_{\rm curv}(\ks)&=\mu\,L_c^2\,\frac{\alpha_2}{4}\, \|\skew \Grad{{\rm curl}\, \f u}\|^2=\mu\,L_c^2\,\alpha_2\, \|\skew \ks\|^2\, , \quad  \f m=2\,\mu\,L_c^2\,\alpha_2\,\skew\,\ks\,.
\end{align}
In that model the non-local force stress tensor $\widetilde{\fg \tau}$ is skew-symmetric as before, but the couple stress tensor $\f m$ is assumed to be completely skew-symmetric as well. Their reasoning, based on a certain restricted understanding of boundary conditions, is critically discussed in Neff et al.~\cite{NeMueGhiMa15}.
 \end{itemize}
\subsection{A variant of the indeterminate couple stress model with symmetric total force stress}
In Ghiba et al.~\cite{GhiNeMaMue15} the isotropic, linear indeterminate couple stress
model has been modified so as to have symmetric total force stress $\widehat{\fg \sigma}$, while retaining the same weak form of the Euler-Lagrange equations. This is possible since the force stress tensor appearing in the balance of forces is only determined up to a self-equilibrated stress-field $\bar{\fg \sigma}$, i.e.
\begin{align}\label{selfstress}
\Div \widetilde{\fg \sigma} + \f f = 0 \, \Leftrightarrow \, \Div (\widetilde{\fg \sigma} + \bar{\fg \sigma}) + \f f = 0 \, , \quad {\text{for any}} \, \bar{\fg \sigma} \,
{\text{with}} \quad \Div \bar{\fg \sigma}=0 \,.
\end{align}
The curvature energy expression of this new model is
\begin{align}
W_{\rm curv}({\rm D}^2_{\f x} \f u) = \mu \, L_c^2 ( \alpha_1 \, \norm{\dev \sym \Curl(\sym \Grad{ \f u})}^2 + \alpha_2 \, \norm{\skew \Curl(\sym \Grad{ \f u})}^2 ) \, .
\end{align}
The strong form of the new model reads
\begin{align}\label{strongform}
\fbox{
\parbox[][2.0cm][c]{12cm}{$ \Div \widehat{\fg \sigma} + \f f = 0 \, , \qquad \widehat{\fg
\sigma} = \fg \sigma + \widehat{\fg \tau}  \in {\rm Sym}(3)  \quad \text{symmetric total force stress} \\
\\
\fg \sigma = 2 \, \mu \, \sym \Grad{ \f u} + \lambda \tr(\Grad{ \f u}) \id \, ,
\qquad \widehat{\fg \tau} = \sym \Curl(\widehat{\f m}) \\
\\
\widehat{\f m} = 2 \, \mu \, L_c^2 \, ( \alpha_1 \, \dev \sym \Curl(\sym \Grad{ \f u})+
\alpha_2 \, \skew \Curl(\sym \Grad{ \f u})) \,. $}}
\end{align}
The total force stress tensor is now $\widehat{\fg \sigma}=\fg \sigma + \widehat{\fg
\tau}$ and the second order couple stress tensor is $\widehat{\f m}$. Similarly as in the indeterminate couple stress theory we have $\tr(\widehat{\f m})=0$. Compared to the
classical indeterminate couple stress theory one can show that $\Div(\widetilde{\fg
\sigma}-\widehat{\fg \sigma})=0$, as claimed. Thus, eq.\eqref{strongform} is a {\sl
couple stress model with symmetric total force stress $\widehat{\fg \sigma}$ and trace
free couple stress tensor $\widehat{\f m}$.} Moreover, the couple stress tensor $\f m$ can be symmetric itself for the possible choice $\alpha_2=0$.
\subsection{Conformal invariance of curvature in favor of the modified couple stress theory}
An infinitesimal conformal mapping \cite{Neff_Jeong_Conformal_ZAMM08,Neff_Jeong_IJSS09}
preserves angles and shapes of infinitesimal figures to the first order. The
inhomogeneity is therefore only a global feature of the mapping and locally no shear or distortional deformation appears. Therefore we require as particular case that a second gradient model based on couple stresses should not ascribe curvature energy to such deformation modes. Put in other words, we will require that
\begin{align}\label{couplestressresponse}
\fbox{
\parbox[][0.5cm][c]{12cm}{
there should not be any couple stress response under deformations of infinitesimal cubes if the cubes are only rigidly rotated and dilated.
}}
\end{align}
In order to prepare the stage for the subsequent development let us introduce a further axiom which is tacitly assumed in classical mechanics. We call it the {\bf axiom of localized response:}
\begin{align}\label{AxiomLocalResponse}
\fbox{
\parbox[][0.5cm][c]{12cm}{
The constitutive equations can be investigated based on the response of the material on the level of the deformation of infinitesimal cubes.
}}
\end{align}
The axiom of localized response, together with requirement \eqref{couplestressresponse} yields that the couple stress tensor $\f m$ should be independent of conformal curvature.\\

\begin{figure}
\centering
\includegraphics[height=45mm]{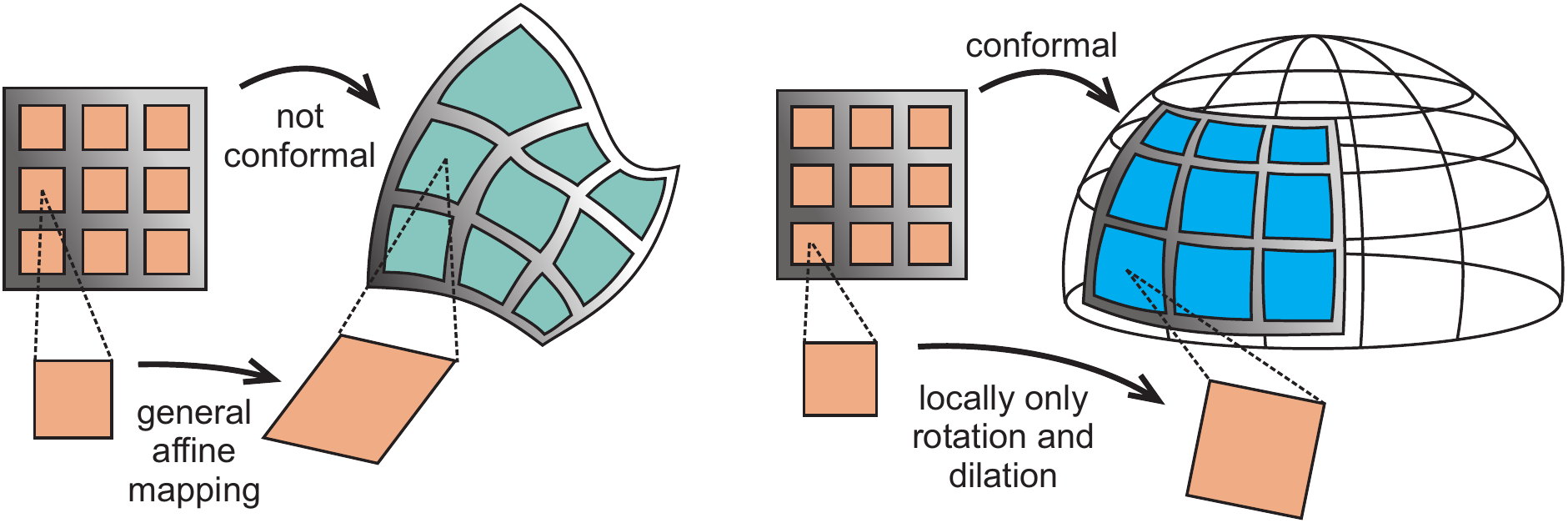}
\caption{General affine (left) and conformal mapping (right), which is locally only rotation and dilation.}\label{PicConvCurv}
\end{figure}

\begin{figure}
\centering
\includegraphics[height=35mm]{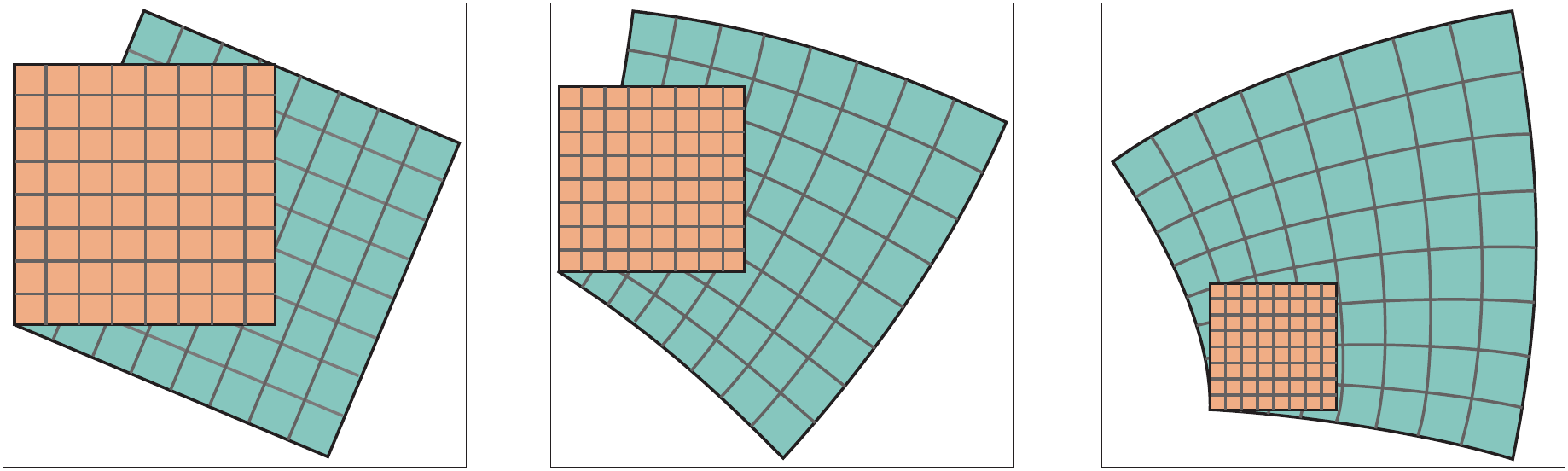}
\caption{Infinitesimal conformal mappings \cite{Neff_Jeong_IJSS09} locally preserve angles and shapes but may be globally  inhomogeneous.}\label{ConformalCurvPic}
\end{figure}

Since this is part of our discussion, we give a short introduction to conformal invariance.
A map $\bar{\f x} = \fg \phi_c(\f x):\mathbb{R}^3\rightarrow\mathbb{R}^3$ is
infinitesimally conformal if and only if its Jacobian satisfies $\nabla \fg \phi_c(\f x)\in
\mathbb{R}\cdot \id+\so(3)$, where $\mathbb{R}\cdot \id+\so(3)$ is the conformal
Lie-algebra. This implies
\cite{Neff_Jeong_Conformal_ZAMM08,Neff_Jeong_IJSS09,Neff_Jeong_bounded_stiffness09} the
representation of that map as a special second order polynomial function
\begin{align}
\fg \phi_c(\f x)=\frac{1}{2}\left(2\langle \axl\widehat{\f W},\f x\rangle \f
x-\axl\widehat{\f W})\|\f x\|^2\right)+[\widehat{p} \, \id+\widehat{\f A}]\cdot \f
x+\widehat{\f b} \, , \qquad \xbar = \fg \phi_c(\f x) \, ,
\end{align}
where $\widehat{\f W},\widehat{\f A}\in \so(3)$, $\widehat{\f b}\in \mathbb{R}^3$,
$\widehat{p}\in \mathbb{R}$ are arbitrary but constant. In Fig.\ref{PicConvCurv} and \ref{ConformalCurvPic} possible deformation modes of $\fg \phi_c$ are drawn. By {\bf conformal invariance} of
the curvature energy, we mean that it vanishes on infinitesimal conformal mapping.%

The axiom of localized response \eqref{AxiomLocalResponse} is tacitly assumed in classical continuum mechanics to avoid the necessity of higher order terms in the Taylor series expansion of the total stress $\widetilde{\fg \sigma}$. The constitutive equations can be discussed (similar to the balance equations as done in section \ref{KapCouples})
based on the response of the material on the level of the deformation of infinitesimal cubes. Similarly, assuming the couple stress tensor $\f m$ to be independent of conformal curvature is such a kind of localization: $\f m(D_x^2\phi_c) = 0 \Rightarrow \f m \in \Sym(3)$.
This is equivalent to
\begin{align}
W_{\rm curv}({\rm D}^2_{\f x} \fg \phi_c)=0 \quad \text{or} \quad \f m({\rm D}^2_{\f x} \fg \phi_c)=0  \qquad \text{for all conformal maps} \, \fg
\phi_c \, .
\end{align}
The classical linear elastic energy still ascribes energy to such a deformation mode but
strictly related only to the bulk modulus
\begin{align}
W_{\rm lin}(\Grad{ \fg \phi_c})=\frac{3\,\lambda+2\,\mu}{2} \, [\tr(\Grad{ \fg \phi_c})]^2 \, ,
\end{align}
i.e. to volumetric deformation parts inherent in $\fg \phi_c$. In case of a classical plasticity formulation with von Mises deviatoric flow rule, conformal mappings are precisely those inhomogeneous mappings that never lead to plastic flow
\cite{Neff_Cosserat_plasticity05} since
\begin{align}
\dev \sym \Grad{\fg \phi_c} \equiv 0 \, .
\end{align}
In that perspective, conformal mappings are ideally elastic in the sense that regardless of how large they are, von Mises plasticity is never triggered. Introducing the displacement field $\f u = \fg \phi(\f x) - \f x :\B \subset \mathbb{R}^3\rightarrow\mathbb{R}^3$, it can be remarked that
\begin{align}
W_{\rm curv} ({\rm D}^2_{\f x} \f u)= W_{\rm curv}(\dev \sym \Grad{ \curl \, \f u})\,, \quad \f m(\ks) =2\,\mu\,L_c^2\,\alpha_1 \, \dev \sym \ks
\end{align}
is conformally invariant. But e.g.~the curvature energy chosen by Hadjesfandiari and Dargush \cite{hadjesfandiari2014evo}
\begin{align}
W_{\rm curv}({\rm D}^2_{\f x} \f u) = W_{\rm curv}(\skew \Grad{ \curl \, \f u})\,,\quad \f m(\ks) =2\,\mu\,L_c^2\,\alpha_2 \, \skew \ks
\end{align}
is not conformally invariant. Our new model (\ref{strongform}) with $\alpha_2=0$ yields
\begin{align}
W_{\rm curv}({\rm D}^2_{\f x}) &= \mu \, L_c^2 \, \alpha_1 \, \norm{\dev \sym \Curl(\sym \Grad{ \f u})}^2 \,, \notag \\
\widehat{\f m} & = 2 \, \mu \, L_c^2 ( \alpha_1 \, \dev \sym \Curl(\sym \Grad{ \f u}) \,,
\end{align}
which is also conformally invariant. Thus, the underlying additional invariance property of the modified couple stress theory is precisely conformal invariance. In the modified couple stress model, these deformations are free of size-effects. Put in other words, the generated couple stress tensor $\widetilde{\f m}$ in the modified couple stress model is zero for this deformation mode, while in the model by Hadjesfandiari and Dargush $\widetilde{\f m}$ is constant and skew-symmetric under infinitesimal conformal mappings.
\section{Another derivation of the equations for the couple stress model}\setcounter{equation}{0}
There are several ways to arrive at the equilibrium equations of the couple
stress model: the formal way postulates energy minimization and results in Euler-Lagrange equation \eqref{CoupleStressBalanceMomentum}. Therefore, constitutive assumptions on the
energy function need to be made. This procedure has been followed in the first part of
this paper. Another route consists in looking at a discrete lattice model, making some assumptions on the next to nearest neighbor interaction and homogenizing the results. This has been followed e.g.~in \cite{Neff_Jeong_IJSS09}. There, the homogenized energy is obtained and equilibrium follows again as an Euler-Lagrange equation.

First, equilibrium equations are obtained by another approach, which is free of constitutive and kinematical assumptions. Thus, it holds for all kind of solid media. We start from a given, generally inhomogeneous total stress distribution $\widetilde{\fg \sigma}$ and postulate equilibrium at subdomains. Subdomains are considered to be infinitesimal and cubic, as traditionally used in classical continuum mechanics. This yields the standard equilibrium equations. Additionally, it can be shown how the equilibrium equations generalize, if the Taylor series expansion of the stress distribution $\widetilde{\fg \sigma}$ allows for higher order terms than usually considered in classical continuum mechanics.

In doing so, we do not introduce other physical quantities besides the total force stress tensor $\widetilde{\fg \sigma}$ for the balance of linear and angular momentum. However, fluctuations of the stress function over infinitesimal cubes are evaluated up to quadratic terms in a spatial Taylor series expansion, which is assumed to be valid within the cubes. On that basis, certain properties of stress can be elaborated with respect to the center of the cube. We extract the couple stress $\f m$ from its assignment to be a stress resultant for the balance of angular momentum, reading
\begin{align}\label{BasisDef_m_from_Sigma}
\f m \colonequals \int_{\partial \B_c} {\rm polar}(\widetilde{\fg \sigma}) \cdot \f n \di A\,.
\end{align}
The polar operator in eq.(\ref{BasisDef_m_from_Sigma}) is not the polar
decomposition but this will be explained later. Our analysis is in principle applicable to any medium, no further constitutive assumptions need to be made. It is perfectly Newtonian in the sense that the whole discussion is based on the statement of balance laws only. However, let us immediately point out the limitations of such an approach:
\begin{itemize}
  \item It is impossible to obtain more general higher gradient models. The interaction   will be limited to some ``rotational" type of effects through the structure of   eq.\eqref{CoupleStressBalanceMomentum}.
  \item It is impossible to obtain a true micromorphic type kinematics since the
  coupling of moment stresses necessary there would also be beyond the presented framework \cite{NeffGhibaMicroModel, GhibaNeffExistence,Neff_Jeong_ZAMP08}.
  \item The approach would offer the possibility for a true Cosserat type kinematics with   independent rotations but does not necessitate these independent degrees of freedom \cite{Neff_ZAMM05}.
\end{itemize}
As a preliminary conclusion we can say: The procedure in the next section \ref{KapTaylor} is one of the many possibilities to motivate the indeterminate couple stress model. Since its assumptions are taken from the traditional balance laws of the kinematics of rigid
bodies, some authors claim that this is ``the one and only" motivation for such a model.
Clearly, we need to dismiss such a strong claim: it is well accepted that continuum
mechanics extends far beyond the kinematics of rigid bodies. On the positive side, this
derivation let us better understand the engineering way of motivation for the couple
stress model. It also allows us to see the fallacy of Yang et al.'s argument later in this paper.

\subsection{Taylor series expansion of total stress}\label{KapTaylor}
We treat the total stress tensor $\widetilde{\fg \sigma}$ as spatial function and no
restrictions on its symmetry apply a priori. We assume that the contact forces acting in
the body are fully described by this total stress $\widetilde{\fg \sigma}$. Cauchy's
principle states that the traction $\fg \sigma_{\rm n}$ on any surface of a body derives
from the force stress $\widetilde{\fg \sigma}$ and the surface normal $\f n$ via
\begin{align}\label{CauchyI}
\fg \sigma_{\rm n} = \widetilde {\fg \sigma} \cdot \f n \,, \qquad \fg \sigma_{\rm n} \in
\R^3 \, , \qquad \widetilde {\fg \sigma} \in \R^{3\times3} \,.
\end{align}
It is a generalization to use this principle also for the couple stress $\f m$ as already proposed by Koiter \cite{Koiter64}, reading
\begin{align}\label{CauchyII}
 \f m_{\rm n} =  \f m \cdot \f n  \,, \qquad \f m_{\rm n} \in \R^3 \, , \qquad
\f m \in \R^{3\times3} \,.
\end{align}
Here, eq.\eqref{CauchyI} and \eqref{CauchyII} are axiomatic in nature. However, the similarity of $\widetilde {\fg \sigma}$ and $\f m$ concerning Cauchy's principle implies that stress and couple stress relate to the same physical quantity: the bonding force between neighboring material points. In our opinion, couple stress represents a certain kind of non-local bonding force of the force stress function, being neglected
within a local continuum formulation. The effect of non-local bonding forces is due to the inhomogeneity of the force-stress function. Thus,  we consider the split of this stress function into several parts subsequently.\\

\begin{figure}
\centering
\includegraphics[height=35mm]{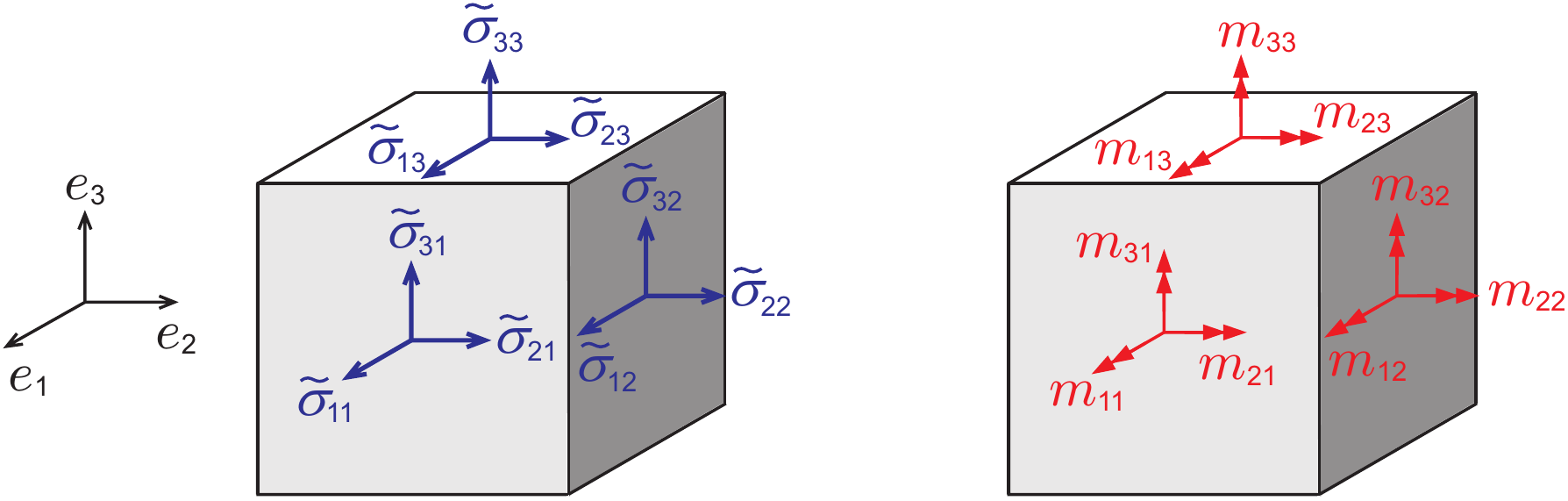}
\caption{Faces of a finite cube $\B_c$ showing components of the force stress tensor $\widetilde{\fg \sigma}$ and the couple stress tensor $\f m$, respectively.}\label{DefSigmaM}
\end{figure}

In accordance with \cite{Truesdell65}, the index $i$ of stress components $\sigma_{ij}$
characterizes the component $(\sigma_n)_i$ of the subsequent force traction in the direction of the associated base vector $\f e_i$, and the index $j$ characterizes the plane that $\fg \sigma_n$ is acting upon\footnote{Some authors reverse this convention by identifying the first index with the plane and the second index with the vector component.}, i.e.~$\f n = \f e_i$ specifies the direction of the plane. The same convention is considered for the couple stress tensor \cite{Toupin64}, see Fig.\ref{DefSigmaM}.

\begin{figure}
\centering
\includegraphics[height=30mm]{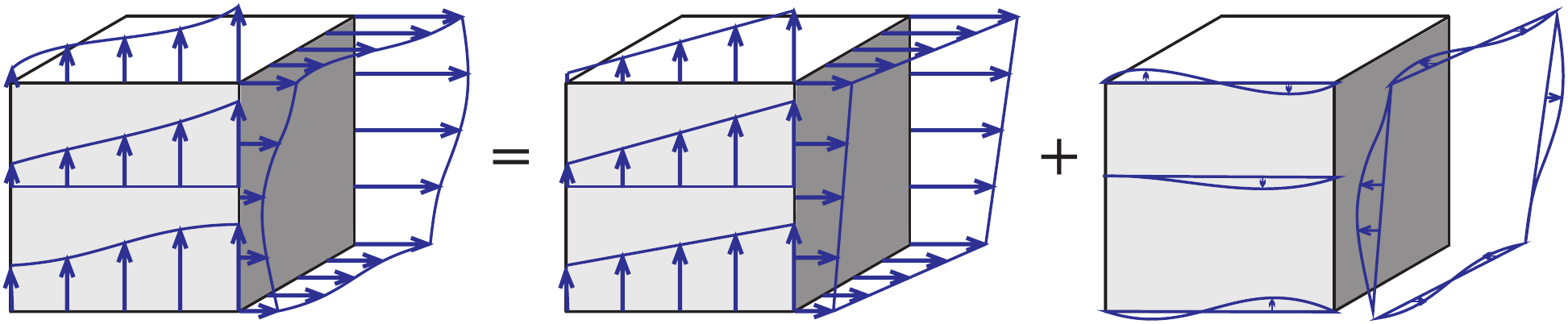}
\caption{Left: Stress fluctuation for a tangential and a normal stress component. Taylor series expansion results in linear (middle) and higher order terms (right).}\label{StressFluctuation}
\end{figure}

The components of the force stress tensor are usually drawn in a simple way at single points centered on each face of Cauchy's cube $\B_c$. However, the total stress tensor depends on its position $\f x$ in space. Generally, fluctuations may appear from point to point, exemplarily sketched for a tangential and a normal stress component in Fig.\ref{StressFluctuation}. Reducing the cube $\B_c$ to the size of a point gives some motivation for the simplified representation with single arrows centered on faces.
\begin{remark}
Even in classical continuum mechanics, stress components generally need to vary
linearly between opposite faces of the infinitesimal cube to appear in the balance
equation of linear momentum, as will be shown below. On the other hand, the same stress
components are treated as constant on faces, where they appear as traction $\widetilde{\fg \sigma} \cdot \f n$. Such a directional selection of stress gradients eliminates moment couples in the stress function, which is due to excluding couple stress in the model.
\end{remark}

\noindent
Next, we discuss this issue in detail with the help of the second-order Taylor series expansion of the total stress
\begin{align}\label{TaylorStress}
\widetilde{\fg \sigma}(\f x_0 + \Delta \f x)
        &= \underbrace{\widetilde{\fg \sigma}(\f x_0)}_{\widetilde{\fg \sigma}^0} +
        \underbrace{\frac{\pt \widetilde{\fg \sigma}(\f x_0)}{\pt x_1}\, \Delta x_1+
        \frac{\pt \widetilde{\fg \sigma}(\f x_0)}{\pt x_2} \, \Delta x_2+
        \frac{\pt \widetilde{\fg \sigma}(\f x_0)}{\pt x_3} \, \Delta x_3}_{\TD . \Deltax}\notag \\
        &+\underbrace{\frac{\pt^2 \widetilde{\fg \sigma}(\f x_0)}{\pt x_1 \, \pt x_2} \, \Delta x_1\, \Delta x_2+
        \frac{\pt^2 \widetilde{\fg \sigma}(\f x_0)}{\pt x_1 \, \pt x_3} \, \Delta x_1 \, \Delta x_3+
        \frac{\pt^2 \widetilde{\fg \sigma}(\f x_0)}{\pt x_2 \, \pt x_3} \, \Delta x_2 \, \Delta x_3}_{\TDb  . \Deltaxb} \notag \\
        &+\underbrace{\12 \frac{\pt^2 \widetilde{\fg \sigma}(\f x_0)}{\pt x_1^2} \, \Delta x_1^2 +
        \12 \frac{\pt^2 \widetilde{\fg \sigma}(\f x_0)}{\pt x_2^2} \, \Delta x_2^2+
        \12 \frac{\pt^2 \widetilde{\fg \sigma}(\f x_0)}{\pt x_3^2} \, \Delta x_3^2}_{\TDq . \Deltaxq}
        + \, o(\Delta \f x^3, \, \Delta \f x^4, \, ...) \,,
\end{align}
with $\Delta \f x = \f x - \f x_0$ describing the distance from the center to any point
of the cube $\B_c$.

The derivatives of the stress function $\TD$, $\TDb$, and $\TDq$ are evaluated in the center of the cube and constant in $\B_c$ and on $\partial \B_c$. As products with $\Deltax$, $\Deltaxb$, and $\Deltaxq$, a fully bilinear representation of stress is given after neglecting higher order terms $o(\Delta \f x^3, \, \Delta \f x^4, \, ...)$. We split second order derivatives in two terms. As we show later, this split is motivated by different effects of each term concerning the balance of angular momentum.\\

\begin{figure}
\centering
\includegraphics[height=30mm]{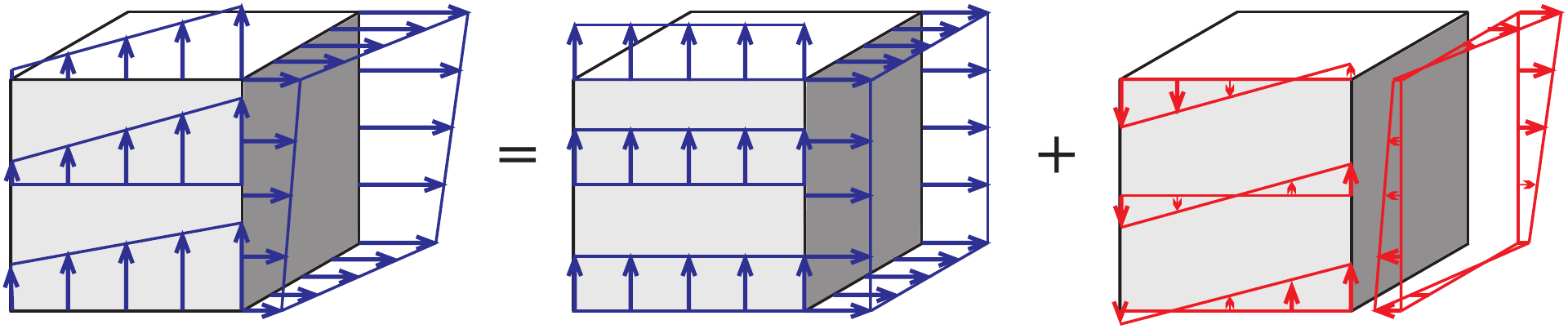}
\caption{Decomposition of linearized stress components (left) in constant (middle) and linear (right) parts.}\label{StressLinearized}
\end{figure}

To illustrate the decomposition of the total stress function  $\widetilde{\fg \sigma}$, a tangential and a normal stress component in Fig.\ref{StressFluctuation}a are exemplarily decomposed into a linear function in Fig.\ref{StressFluctuation}b and higher order terms in Fig.\ref{StressFluctuation}c.
Even linear fluctuations of stress, represented by $\TD . \Deltax$, generally yield couples regarding the center of faces. In Fig.\ref{StressLinearized} a linearized tangential and normal stress component are drawn. Their decomposition into constant and linear terms are shown. Obviously, constant stress does not generate a couple concerning the center of the face, where it acts on. This is in contrast to the linear fluctuation, which obviously results in a couple. Since this physical effect of generating a couple or not is essential for this paper, we use the following terminology:
\begin{align}
\fbox{
\parbox[][0.5cm][c]{12cm}{
Terms of stresses generating tractions $\widetilde{\fg \sigma} \cdot \f n_i$ such
that a couple emerges with respect to the center of the surface $i$ are called {\bf polar}.}}
\end{align}
Let us repeat that in classical continuum mechanics, stress components are treated as linear functions between opposite faces of an infinitesimal cube. But the samestress components are assumed to be constant on faces, where they appear as traction $\widetilde{\fg \sigma} \cdot \f n$. Therefore, by completely keeping such linear terms ofstress, exemplarily shown in Fig.\ref{StressLinearized}, the classical approach of
continuum mechanics is extended. We will show in the next section that certain terms of the Taylor series expansion yield contributions to the balance of linear momentum. Similarly, some terms yield contributions to the balance of angular momentum. This is why we define a second terminology:
\begin{align}
\fbox{
\parbox[][0.5cm][c]{12cm}{
Terms of stresses generating neither polar tractions nor contributing to the
balance of angular momentum are called {\bf nonpolar}.}}
\end{align}
It is crucial for this paper that linear and bilinear terms from the Taylor series expansion in eq.\eqref{TaylorStress} will be considered in order to find the origins and properties of couple stress.
\subsection{Discussion of origins and properties of the couple stress}\label{KapCouples}
\begin{figure}
\centering
\includegraphics[height=45mm]{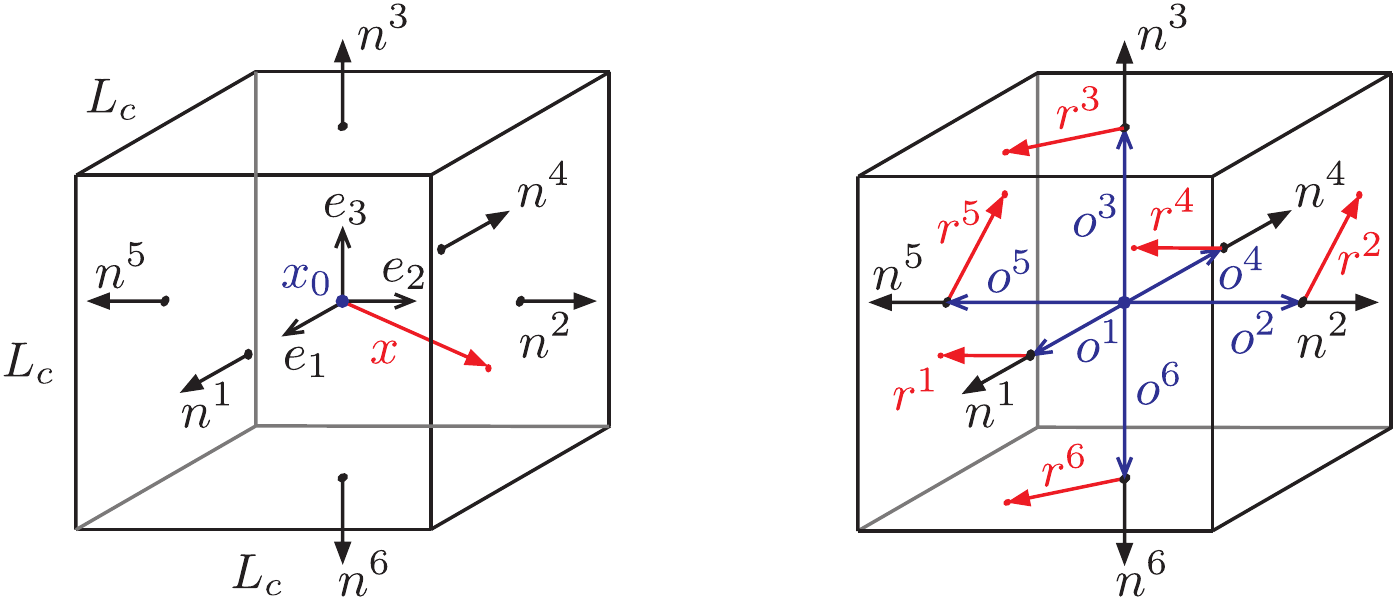}
\caption{Finite cube $\B_c$ with cartesian coordinate system $\f e_i$, center $x_0$, dimension $L_c$, and  normal vectors $\f n^I$ on faces $\partial \B_c^I$. Tangential vectors $\f r^I$ are sketched via offset $\f o^I$ onto the corresponding surface $\partial B_c^I$.}\label{CubeKOS}
\end{figure}

From now on, a local cartesian coordinate system with basis vectors $\f e_i$ aligned to
the edges of the finite cube $\B_c$ is used. The origin $\f x_0$ is considered in the
center of the cube, as shown in Fig.\ref{CubeKOS}. Thus, the increments $\Delta x_1$,
$\Delta x_2$, and $\Delta x_3$ in eq.\eqref{TaylorStress} are aligned along the cartesian coordinates $x_1$, $x_2$, $x_3$. Dimensions of $\B_c$ are limited to the finite length $L_c$ such that $\Delta x_i = x_i \in [-L_c /2 \, , \, L_c/2] \, , i=1,2,3$. The volume of the cube is given by $V_c =L_c^3$. The six faces of the cube are indicated according to their normal vectors:
\begin{align}\label{faceparametrization1}
\begin{array}{lrllrllrl}
  \f n_1 = & \f e_1 & \text{on} \, \B_c^1 \, , \quad & \f n_2 = & \f e_2 & \text{on} \, \B_c^2 \, , \quad  &  \f n_3 = & \f e_3 & \text{on} \, \B_c^3 \, ,  \\
  \f n_4 = & -\f e_1 & \text{on} \, \B_c^4 \, , \quad & \f n_5 = & - \f e_2 & \text{on} \, \B_c^5 \, , \quad  &  \f n_6 = &
  -\f e_6 & \text{on} \, \B_c^6 \, .
\end{array}
\end{align}
The cubes faces are parameterized by cartesian coordinates defining the tangent vectors
\begin{align}\label{faceparametrization}
\f r_1= \f r_4 = \left(
                   \begin{array}{c}
                     0 \\
                     x_2 \\
                     x_3 \\
                   \end{array}
                 \right) \, , \quad \quad
\f r_2= \f r_5 = \left(
                   \begin{array}{c}
                     x_1 \\
                     0 \\
                     x_3 \\
                   \end{array}
                 \right) \, , \quad \quad
\f r_3= \f r_6 = \left(
                   \begin{array}{c}
                     x_1 \\
                     x_2 \\
                     0 \\
                   \end{array}
                 \right) \, .
\end{align}
Each parametrization $\f r_i$ is face centered and perpendicular to the
normal $\f n_i$ and hence tangent to the surface itself.
First, the {\rm balance of linear momentum} is discussed in view of using the Taylor series expansion of eq.(\ref{TaylorStress}). By Gauss's theorem, the sum of tractions on all faces of the cube reads
\begin{align}\label{SigmaBalanceMomentum}
\sum_{i=1}^6 \int_{\partial \B^i_c} \widetilde{\fg \sigma}(\f x_0 + \Delta \f x) \cdot \f n_i \di A = \int_{ \B_c} \Div \widetilde{\fg \sigma}(\f x_0 + \Delta \f x) \, \di V \,.
\end{align}
Taking the divergence of the second-order Taylor series expansion of $\widetilde{\fg \sigma}$ from eq.\eqref{TaylorStress} gives
\begin{align}\label{DivTaylor}
\Div & \widetilde{\fg \sigma}(\f x_0 + \Delta \f x) = [ \widetilde{\sigma}_{i1,1}(\f x_0) + \widetilde{\sigma}_{i2,2}(\f x_0) + \widetilde{\sigma}_{i3,3}(\f x_0) + (
\widetilde{\sigma}_{i1,11}(\f x_0) + \widetilde{\sigma}_{i2,12}(\f x_0) +
\widetilde{\sigma}_{i3,13}(\f x_0)) \Delta x_1
\notag \\
&+ ( \widetilde{\sigma}_{i2,22}(\f x_0) + \widetilde{\sigma}_{i1,12}(\f x_0)
+\widetilde{\sigma}_{i3,23}(\f x_0)) \Delta x_2 + ( \widetilde{\sigma}_{i3,33}(\f x_0) +
\widetilde{\sigma}_{i1,13}(\f x_0) +\widetilde{\sigma}_{i2,23}(\f x_0)) \Delta x_3 ] \,
\f e_i \,.
\end{align}
Naturally, constant terms $\widetilde{\fg \sigma}^0$ do not appear in eq.\eqref{DivTaylor}, which, by recalling again Gauss' theorem, is in accordance with
\begin{align}\label{BOM_sigma0}
\sum_{i=1}^6 \int_{\partial \B^i_c} \widetilde{\fg \sigma}^0 \cdot \f n_i \di A =0 \,.
\end{align}
The last eq.\eqref{BOM_sigma0} implies that $\widetilde{\fg \sigma}^0$ yields opposite constant tractions at opposite faces, due to opposing normal vectors defined in \eqref{faceparametrization1}. Of course, it is a well known fact that $\widetilde{\fg \sigma}^0$ has no influence on the balance of linear momentum, which only involve incremental quantities.\\

Linear terms concerning one component of $\Delta \f x$ appear in eq.\eqref{DivTaylor}. Thus, symmetric bounds of integration $- \di x/ 2$ to $\di x/ 2$ cancel each other out when performing the body volume integration in eq.\eqref{SigmaBalanceMomentum}. Consequently, the terms $\TDb . \Deltaxb$ and $\TDq . \Deltaxq$ do not contribute to the balance of linear momentum. In formulas:
\begin{align}\label{BOM_biquad}
\sum_{i=1}^6 \int_{\partial \B^i_c} & (\TDb . \Deltaxb + \TDq . \Deltaxq) \cdot \f n_i \di A \notag \\
&= \int_{ \B_c} \Div (\TDb . \Deltaxb + \TDq . \Deltaxq) \, \di V = 0 \,.
\end{align}
Thus, we can further work on eq.\eqref{SigmaBalanceMomentum} which, using  \eqref{BOM_sigma0} and \eqref{BOM_biquad} implies
\begin{align}\label{SigmaBalanceMomSpec}
\sum_{i=1}^6 \int_{\partial \B^i_c} & \widetilde{\fg \sigma}(\f x_0 + \Delta \f x) \cdot \f n_i \di A \notag \\
&= \sum_{i=1}^6 \int_{\partial \B^i_c} (\widetilde{\fg \sigma}^0 + \TD.\Deltax
+ \TDb.\Deltaxb +  \TDq.\Deltaxq) \cdot \f n_i \di A \notag\\
&= \int_{\B_c} \Div (\TD.\Deltax) \di V \,,
\end{align}
or equivalently
\begin{align}\label{SigmaBalanceMomSpecX}
\sum_{i=1}^6 \int_{\partial \B^i_c} & \widetilde{\fg \sigma}(\f x_0 + \Delta \f x) \cdot \f n_i \di A = \int_{\B_c} \Div (\TD.\Deltax) \di V  \notag\\
& =  \int_{\B_c} \underbrace{[\widetilde{\sigma}_{i1,1}(\f x_0)+\widetilde{\sigma}_{i2,2}(\f x_0)+\widetilde{\sigma}_{i3,3}(\f x_0)]}_{\rm const.} \, \f e_i \di V = V_c \,\, \Div(\widetilde{\fg \sigma}(\f x_0)) \,,
\end{align}
where $V_c$ is the volume of the cube. The result in eq.\eqref{SigmaBalanceMomSpecX} shows that even for the bilinear approximation of stress within the Cauchy cube $\B_c$, the evaluation of $\Div(\widetilde{\fg \sigma})$ in its center $\f x_0$ is sufficient to represent the sum of tractions on $\partial \B_c$\footnote{However, including higher order terms $o(\Delta \f x^3)$ would let appear higher order derivatives of $\widetilde{\sigma}$.}. The reader should be aware of the fact that we do not postulate
$\Div(\widetilde{\fg \sigma}(\f x)) = \Div(\widetilde{\fg \sigma}(\f x_0))$ in $\B_c$ and, therefore, the result in eq.(\ref{SigmaBalanceMomSpecX}) is not trivial and only valid if we limit the Taylor series expansion in eq.\eqref{TaylorStress} up to second order.\\

Further, the balance of angular momentum contains only some terms of $\TD$, being constant in eq.\eqref{SigmaBalanceMomSpecX} and reading $[\widetilde{\sigma}_{i1,1}(\f x_0)+\widetilde{\sigma}_{i2,2}(\f x_0)+\widetilde{\sigma}_{i3,3}(\f x_0)]$. Thus, the argument $\TD.\Deltax$ in eq.\eqref{SigmaBalanceMomSpec} can be decomposed into
\begin{align}\label{SplitSigma}
\TD.\Deltax = \underbrace{\TDnp . \Deltax}_{\rm nonpolar} + \underbrace{\TDp .
\Deltax}_{\rm polar} \, ,
\end{align}
with
\begin{align}\label{LinStressNonpolar}
\TDnp . \Deltax\colonequals\left( \begin{array}{ccc}
   \widetilde{\sigma}_{11,1} & 0 & 0 \\
   \widetilde{\sigma}_{21,1} & 0 & 0 \\
   \widetilde{\sigma}_{31,1} & 0 & 0 \\
   \end{array}
 \right) \Delta x_1
   +\left( \begin{array}{ccc}
   0 & \widetilde{\sigma}_{12,2} & 0 \\
   0 & \widetilde{\sigma}_{22,2} & 0 \\
   0 & \widetilde{\sigma}_{32,2} & 0 \\
   \end{array}
 \right) \Delta x_2
   +\left( \begin{array}{ccc}
   0 & 0 & \widetilde{\sigma}_{13,3} \\
   0 & 0 & \widetilde{\sigma}_{23,3} \\
   0 & 0 & \widetilde{\sigma}_{33,3} \\
   \end{array}
 \right) \Delta x_3 \, ,
\end{align}
and
\begin{align}\label{LinStressPolar}
\TDp . \Deltax \colonequals \!\!
   \left( \begin{array}{ccc}
     0 & \widetilde{\sigma}_{12,1} & \widetilde{\sigma}_{13,1} \\
     0 & \widetilde{\sigma}_{22,1} & \widetilde{\sigma}_{23,1} \\
     0 & \widetilde{\sigma}_{32,1} & \widetilde{\sigma}_{33,1} \\
   \end{array}
 \right) \Delta x_1
   +\! \left( \begin{array}{ccc}
    \widetilde{\sigma}_{11,2} & 0 & \widetilde{\sigma}_{13,2} \\
    \widetilde{\sigma}_{21,2} & 0 & \widetilde{\sigma}_{23,2} \\
    \widetilde{\sigma}_{31,2} & 0 & \widetilde{\sigma}_{33,2} \\
   \end{array}
 \right) \Delta x_2
   +\! \left( \begin{array}{ccc}
   \widetilde{\sigma}_{11,3} & \widetilde{\sigma}_{12,3} & 0 \\
   \widetilde{\sigma}_{21,3} & \widetilde{\sigma}_{22,3} & 0 \\
   \widetilde{\sigma}_{31,3} & \widetilde{\sigma}_{32,3} & 0 \\
   \end{array}
 \right) \Delta x_3 \,.
\end{align}
Since $\TDnp . \Deltax$ solely affects the balance of linear momentum, we can simplify eq.\eqref{SigmaBalanceMomentum} with help of our decomposition in eq.\eqref{SplitSigma} reading
\begin{align}\label{SigmaBalanceMomSpec2}
\sum_{i=1}^6 \int_{\partial \B^i_c} \widetilde{\fg \sigma}(\f x_0 + \Delta \f x) \cdot \f n_i \, \di A = \int_{\B_c} \Div (\TD.\Deltax) \, \di V = \int_{\B_c} \Div (\TDnp.\Deltax)
\, \di V \,.
\end{align}
The index ``np" and ``p" abbreviates ``nonpolar" and ``polar", respectively, which is
anticipated from subsequent results. According to eq.\eqref{SigmaBalanceMomSpec2} and
considering a constant net force $\f f$ within $\B_c$, the balance of linear momentum is
given by\footnote{Here, we neglect dynamical effects.}
\begin{align}\label{Impulsbilanz1}
\sum_{i=1}^6 \int_{\partial \B^i_c} \widetilde{\fg \sigma}(\f x_0  + \Delta \f x) \cdot \f n_i \, \di A & + \int_{ \B_c} \f f \, \di V = 0 \, \Leftrightarrow \, \int_{ \B_c} [ \Div \TDnp . \Deltax + \f f ] \, \di V = 0 \,\notag\\
& \Leftrightarrow \, [ \Div (\widetilde{\fg \sigma}(\f x_0)) + \f f ] \, V_c = 0 \,.
\end{align}
The well known result in eq.\eqref{Impulsbilanz1} shows that our barycentric coordinate system is appropriate to cover classical results. Since only terms of $\TDnp$ contribute to the balance of linear momentum, it legitimates the classical continuum theory to handle stress components constant on faces, where they appear as traction $\widetilde{\fg \sigma} \cdot \f n$. The term $\TDnp.\Deltax$ is sufficient to obtain the classical equation for the balance of linear momentum.\\

For the following discussion, we need to specify and extend our terminology from section \ref{KapTaylor}, which was motivated by two physical effects:
\begin{itemize}
\item Traction from total stress may generate a couple concerning the mid point of a corresponding face $\partial \B_c^i$.
\item Traction from total stress may generate a couple concerning the mid point of the cube $\B_c$, which affects the balance of angular momentum of $\B_c$.
\end{itemize}
Combining both effects leads to altogether four cases. Terms of stresses in our expansion can be distinguished by four cases of polarity:
\begin{align}\label{PolarDef}
\fbox{
\parbox[][5cm][c]{14cm}{
\begin{enumerate}
\item{{\bf polar:}} Terms of total stress generating tractions $\widetilde{\fg \sigma} \cdot \f n_i$ such that couples emerge with respect to $\f r_i$ on $\partial \B_c^i$ are called polar.
\item{{\bf nonpolar:}} Terms of total stress which are not polar by definition 1 nor contributing to the balance of angular momentum with respect to the center of a cube $\B_c$ are called nonpolar.
\item{{\bf semipolar:}} Terms of total stress which are not polar by definition 1 but contributing to the balance of angular momentum with respect to the center of a cube $\B_c$ are called semipolar.
\item{{\bf bipolar:}} Terms of total stress which are polar by definition 1 and contribute to the balance of angular momentum with respect to the center of a cube $\B_c$ are called bipolar.
\end{enumerate}
}}
\end{align}
Naturally, the center of faces and the center of the cube are neutral points of rotation to define couples on faces and the balance of angular momentum. Thus, the lever arms of tractions acting on the cubes faces $\partial \B^i_c$ are given by $\f r_i$ and $\f x_P$, respectively. Face centered lever arms $\f r_i$ and the origin vector $\f x_P$ are sketched in Fig.\ref{CubeKOS}.\\

We start our discussion with the polarity properties of $\widetilde{\fg \sigma}^0$. Since $\widetilde{\fg \sigma}^0  \cdot \f n_i$ is constant on $\partial \B_c$ it yields the resulting couple on each face $\partial \B_c^i$ with respect to their center to be zero, reading
\begin{align}\label{Sigma0Nonpolar}
\int_{\partial \B^i_c}  \f r_i \times \widetilde{\fg \sigma}^0  \cdot \f n_i \di A = 0 \, , \quad \quad i=1,...,6 \,.
\end{align}
Thus $\widetilde{\fg \sigma}^0$ is not polar according to our Definition \eqref{PolarDef}$_1$. On the other hand, it is not nonpolar but semipolar since skew-symmetric parts of $\widetilde{\fg \sigma}^0$ contribute to the balance of angular momentum via
\begin{align}\label{Sigma0Semipolar}
\sum_{i=1}^6 \int_{\partial \B^i_c}  \f x_P \times \widetilde{\fg \sigma}^0  \cdot \f n_i \di A = \left(
  \begin{array}{c}
    \widetilde{\sigma}^0_{32}- \widetilde{\sigma}^0_{23} \\
    \widetilde{\sigma}^0_{13}- \widetilde{\sigma}^0_{31} \\
    \widetilde{\sigma}^0_{21}- \widetilde{\sigma}^0_{12} \\
  \end{array}
\right) \, L_c^3 = - \fg \epsilon : \widetilde{\fg \sigma}^0 \, V_c = \fg \epsilon :
\fg (\widetilde{\fg \sigma}^0)^T = 2 \, \axl \skew  \widetilde{\fg \sigma}^0  \, V_c \,.
\end{align}
\begin{remark} In accordance with our terminology on polarity in Definition (\ref{PolarDef}) we conclude that constant stress is in general semipolar. In case of $\widetilde{\fg \sigma}\in \Sym(3)$, constant stress is indeed nonpolar.\end{remark}

Next, using the split from eq.\eqref{SplitSigma} to analyze stress gradients from $\TDnp$ reveals
\begin{align}\label{Lin_d_Nonpolar}
\int_{\partial \B^i_c} \f r_i \times (\TDnp . \Deltax) \cdot \f n_i \di A = 0 \,,\quad \quad i=1,...,6 \,,
\end{align}
and
\begin{align}\label{Lin_d_NonSemipolar}
\sum_i \int_{\partial \B^i_c} \f x_P \times (\TDnp . \Deltax)  \cdot \f n_i \di A = 0 \,.
\end{align}
Since the surface traction in eq.\eqref{Lin_d_Nonpolar} and \eqref{Lin_d_NonSemipolar} has generally no polar effect we conclude: \begin{remark} Stress gradients specified by $\TDnp . \Deltax$ from eq.\eqref{LinStressNonpolar} are generally nonpolar even if the total stress $\widetilde{\fg \sigma} \notin \Sym(3)$.
\end{remark}

In contrast, $\TDp . \Deltax$ given by eq.\eqref{LinStressPolar} yields
\begin{align}\label{Lin_p_polar}
\int_{\partial \B^1_c}  \f r_1 \times (\TDp . \Deltax)  \cdot \f n_1 \di A = \frac{L_c^4}{12}\left(
  \begin{array}{c}
    \widetilde{\sigma}_{31,2}- \widetilde{\sigma}_{21,3} \\
    \widetilde{\sigma}_{11,3} \\
    - \widetilde{\sigma}_{11,2} \\
  \end{array}
\right)\, , \notag \\
\int_{\partial \B^2_c}  \f r_2 \times (\TDp . \Deltax)  \cdot \f n_2 \di A = \frac{L_c^4}{12}\left(
  \begin{array}{c}
    - \widetilde{\sigma}_{22,3} \\
    \widetilde{\sigma}_{12,3}- \widetilde{\sigma}_{32,1} \\
    \widetilde{\sigma}_{22,1} \\
  \end{array}
\right)\, ,\notag \\
\int_{\partial \B^3_c}  \f r_3 \times (\TDp . \Deltax)  \cdot \f n_3 \di A = \frac{L_c^4}{12}\left(
  \begin{array}{c}
    \widetilde{\sigma}_{33,2} \\
    - \widetilde{\sigma}_{33,1} \\
    \widetilde{\sigma}_{23,1}- \widetilde{\sigma}_{13,2} \\
  \end{array}
\right)\, , \\
\int_{\partial \B^4_c}  \f r_4 \times (\TDp . \Deltax)  \cdot \f n_4 \di A = -
\int_{\partial \B^1_c}  \f r_1 \times (\TDp . \Deltax)  \cdot \f n_1 \di A \, , \notag \\
\int_{\partial \B^5_c}  \f r_5 \times (\TDp . \Deltax)  \cdot \f n_5 \di A = -
\int_{\partial \B^2_c}  \f r_2 \times (\TDp . \Deltax)  \cdot \f n_2 \di A \, , \notag \\
\int_{\partial \B^6_c} \f r_6 \times (\TDp . \Deltax)  \cdot \f n_6 \di A = - \int_{\partial
\B^3_c}  \f r_3 \times (\TDp . \Deltax)  \cdot \f n_3 \di A_3 \,. \notag
\end{align}
Thus, stress gradients $\TDp.\Deltax$ are polar but do not contribute to the balance of
angular momentum. This results from
\begin{align}\label{Lin_p_notBipolar}
\sum_{i=1}^6 \int_{\partial \B^i_c}  \f x_P \times (\TDp . \Deltax)  \cdot \f n_i
\di A = 0 \,,
\end{align}
which in turn follows from $\f n_1 = - \f n_4$, $\f n_2 = - \f n_5$, and $\f n_3= - \f n_6$. Such a behavior is similar to constant internal stress $\widetilde{\fg \sigma}_0$, which does not influence linear momentum. Therefore we may claim that the physical quantities detected in eq.\eqref{Lin_p_polar} are the constant components of couple stress, reading
\begin{align}\label{Def_m_from_Sigma}
\fbox{
\parbox[][1.0cm][c]{10cm}{$
\f m \colonequals \displaystyle \frac{L_c^2}{12} \left(
  \begin{array}{ccc}
    ( \widetilde{\sigma}_{31,2}- \widetilde{\sigma}_{21,3}) & - \widetilde{\sigma}_{22,3} & \widetilde{\sigma}_{33,2} \\
    \widetilde{\sigma}_{11,3} & ( \widetilde{\sigma}_{12,3}- \widetilde{\sigma}_{32,1}) & - \widetilde{\sigma}_{33,1} \\
    - \widetilde{\sigma}_{11,2} & \widetilde{\sigma}_{22,1} & ( \widetilde{\sigma}_{23,1}- \widetilde{\sigma}_{13,2})
  \end{array}
\right) \, . $}}
\end{align}
Again, couple stress components found from stress gradients $\TDp.\Deltax$ are constant
within $\B_c$ and their convention of sign is as illustrated in Fig.\ref{DefSigmaM}.
Since constant couple stress should not contribute to the balance of angular momentum, eq.\eqref{Lin_p_notBipolar} is in accordance with physical requirements.\\

\begin{figure}
\centering
\includegraphics[height=30mm]{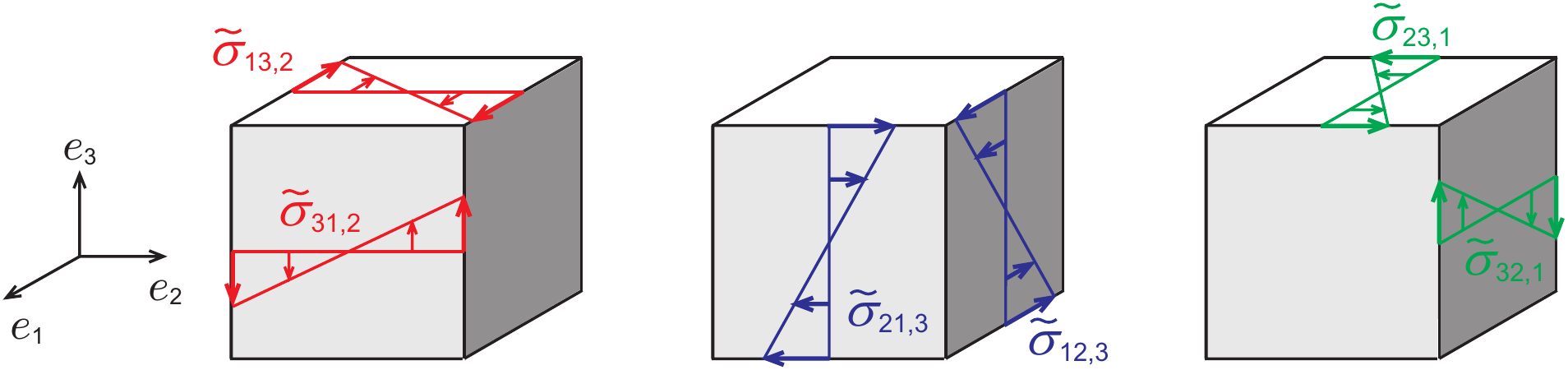}
\caption{Varying shear stress components yielding main diagonal components of couple stress. For a symmetric total stress tensor the sketched stress gradients with corresponding color need to be equal.}
 \label{M_DiagCompo}
\end{figure}

Next, let us discuss the formula for the couple stress tensor $\f m$ from eq.\eqref{Def_m_from_Sigma} in detail. The main diagonal components in eq.\eqref{Def_m_from_Sigma} represent couple stress normal to $\B_c$ according to eq.\eqref{CauchyII}. Thus, couple stress normal to $\B_c$ is caused by fluctuation of shear components in $\widetilde{\fg \sigma}$ as illustrated in Fig.\ref{M_DiagCompo}. It is an interesting result of this derivation that symmetric total force stress  $\widetilde{\fg \sigma}$ yields trace free couple stress $\f m$ when it takes the form given in eq.\eqref{Def_m_from_Sigma}. In fact:
\begin{align}\label{tr_m_from_Sigma}
 \widetilde{\fg \sigma} \in \Sym(3)
 \overset{\text{\eqref{Def_m_from_Sigma}}}{\Longrightarrow}
  \tr (\f m) &= \frac{L_c^2}{12}
    ( \widetilde{\sigma}_{31,2}- \widetilde{\sigma}_{13,2}+
    \widetilde{\sigma}_{12,3}- \widetilde{\sigma}_{21,3}+
    \widetilde{\sigma}_{23,1}- \widetilde{\sigma}_{32,1}) \notag \\
    &= \frac{L_c^2}{12}
    ( (\widetilde{\sigma}_{31}- \widetilde{\sigma}_{13})_{,2}+
    (\widetilde{\sigma}_{12}- \widetilde{\sigma}_{21})_{,3}+
    (\widetilde{\sigma}_{23}- \widetilde{\sigma}_{32})_{,1})
    \equiv 0 \, .
\end{align}
\begin{remark} Since the main diagonal components of the couple stress tensor in eq.\eqref{Def_m_from_Sigma} are associated to space variations of shear stress components, they are not connected to any volue change of an elastic body. The corresponding deformation modes are pure twist of the finite cube.
\end{remark}
In Section \ref{KapTraceFree} we present a symmetric force stress function which yields
trace free and constant couple stress $\f m$. It is another interesting result of our
derivation that:
\begin{remark} The indeterminate couple stress model must have a trace free couple stress tensor $\tr(\f m) = 0$ provided that one assumes a symmetric total stress tensor $\widetilde{\fg \sigma} \in \Sym(3) $ at the beginning. This result is otherwise independent of any further constitutive assumption.
\end{remark}

\begin{figure}
\centering
\includegraphics[height=55mm]{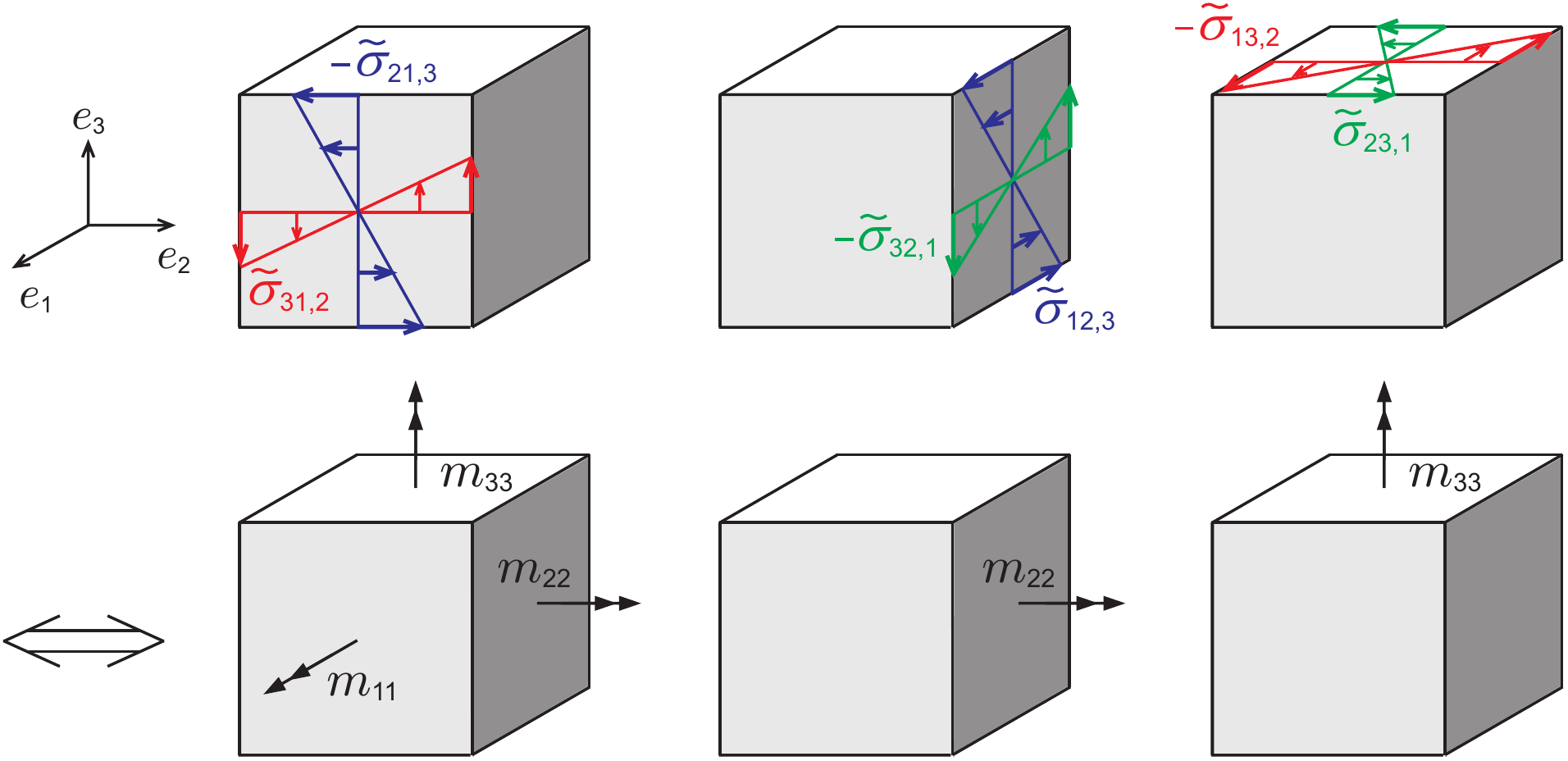}
\caption{Varying shear stress components as origin of main diagonal components of couple stress.}
 \label{MTorsion}
\end{figure}

\begin{figure}
\centering
\includegraphics[height=25mm]{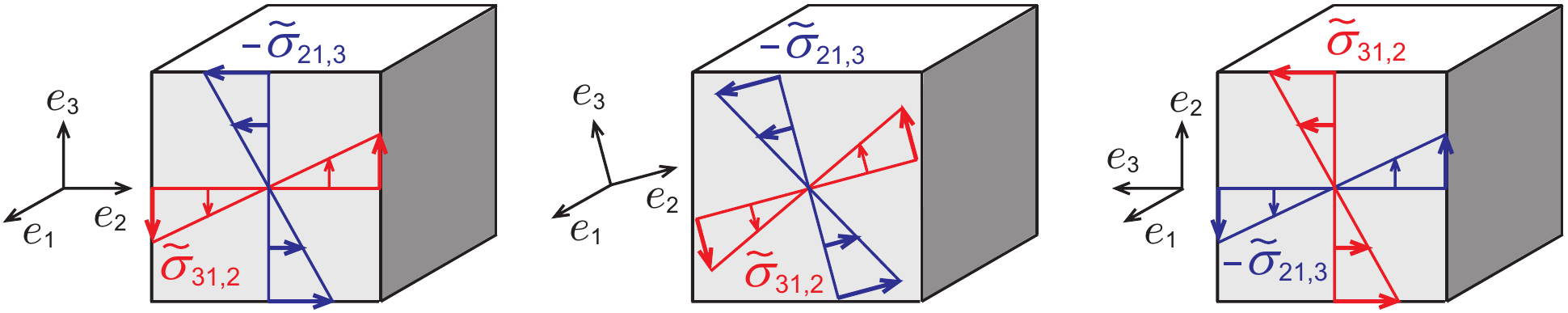}
\caption{Invariance of coplanar stress gradients under rotations around the $\f e_1$-axis.}
 \label{mTorsionInv}
\end{figure}

The off-diagonal components in eq.\eqref{Def_m_from_Sigma} are couples tangential
to $\B_c$ resulting from gradients of normal force stress on $\B_c$. Let us seperate the diagonal (torsion) and off-diagonal (bending) parts
\begin{align}
\f m = \f m^{\rm torsion} + \f m^{\rm bending}
\end{align}
in order to distinguish normal and tangential couples on $\partial \B_c$:
\begin{align}
&\f m^{\rm torsion} \colonequals \frac{L_c^2}{12} \left(
  \begin{array}{ccc}
    ( \widetilde{\sigma}_{31,2}- \widetilde{\sigma}_{21,3}) & 0 & 0 \\
      0 & ( \widetilde{\sigma}_{12,3}- \widetilde{\sigma}_{32,1}) & 0 \\
      0 & 0 & ( \widetilde{\sigma}_{23,1}- \widetilde{\sigma}_{13,2})
  \end{array}
\right) \, , \label{Def_m_torsion} \\
&\f m^{\rm bending} \colonequals \frac{L_c^2}{12} \left(
  \begin{array}{ccc}
    ( 0 & - \widetilde{\sigma}_{22,3} & \widetilde{\sigma}_{33,2} \\
    \widetilde{\sigma}_{11,3} & 0 & - \widetilde{\sigma}_{33,1} \\
    - \widetilde{\sigma}_{11,2} & \widetilde{\sigma}_{22,1} & 0
  \end{array}
\right) \, . \label{Def_m_bending}
\end{align}
In Fig.\ref{MTorsion} the normal couples from eq.\eqref{Def_m_torsion} are drawn. Since
these torsional components have two causes from coplanar stress gradients, let us discuss next if, and how, both are related.
For example, the distribution of shear stresses on $\partial \B^1_c$ is invariant under rotations around the $\f e_1$ axis if $\sigma_{31,2} = - \sigma_{21,3}$, as sketched in Fig.\ref{mTorsionInv}.
With the same invariance argument on $\partial \B^2_c$ and $\partial \B^3_c$, we obtain
truly spatial conditions on gradients of total stress $\widetilde{\fg \sigma}$, reading
\begin{align}\label{mTorsionCriterion}
\sigma_{31,2} = - \sigma_{21,3} \, , \qquad \sigma_{12,3} = - \sigma_{32,1} \, , \qquad
\sigma_{23,1} = - \sigma_{13,2} \, .
\end{align}
In other words, assuming that $\f m^{\rm torsion}$ is given, e.g.~by a balance equation,
the statements of eq.\eqref{mTorsionCriterion} yield \textit{rotationally-invariant coplanar stress gradients} on the Cauchy cube.\\

We proceed our discussion by postulating that
\begin{align}\label{Grad3sph}
\fbox{
\parbox[][1.0cm][c]{14cm}{
in an elastic solid, couple stress $\f m$ should {\sl neither} be connected to the spherical part of the total stress field $\widetilde{\fg \sigma}$ {\sl nor} to its gradient:
\begin{center}
$
\Grad {3 \, \text{sph} \, \widetilde{\fg \sigma}}= \Grad{\tr \widetilde{\fg \sigma}} \nsim f(\f m) \,.
$
\end{center}
}}
\end{align}
This is motivated by two aspects:
\begin{itemize}
  \item {The gradient of dilational stress $\fg \sigma^d = p \, \id$ typically affects the dynamics of fluids where couple stresses are meaningless.}
  \item Conformal mappings, which are generally non-isochoric, should yield no couple stress $\f m$.
\end{itemize}
Thus, the gradient of spherical stress should not be related to $\f m$. However,
\begin{align}\label{Gradtrsigma}
\text{Grad} [\tr \widetilde{\fg \sigma}] = \left(\begin{array}{c} \widetilde{\sigma}_{11,1} \\ \widetilde{\sigma}_{22,2} \\ \widetilde{\sigma}_{33,3} \end{array}\right) + \underbrace{\left(\begin{array}{c} (\widetilde{\sigma}_{22} + \widetilde{\sigma}_{33})_{,1} \\ (\widetilde{\sigma}_{11} + \widetilde{\sigma}_{33})_{,2} \\(\widetilde{\sigma}_{11} + \widetilde{\sigma}_{22})_{,3} \end{array}\right)}_{\displaystyle = \frac{24}{L_c^2} \, \axl(\skew \f m)} = \Div (\text{Diag} \, \widetilde{\fg \sigma}) + \frac{24}{L_c^2} \, \axl(\skew \f m) \,,
\end{align}
contains the skew-symmetric parts of $\f m$. Thus, postulating eq.\eqref{Grad3sph} together with eq.\eqref{Gradtrsigma} necessitates $\f m \in \Sym(3)$, constraining the bending part $\f m^{\text{bending}}$ by three further conditions:
\begin{align}\label{sym_m_from_Sigma}
 &\widetilde{\sigma}_{11,3} = - \widetilde{\sigma}_{22,3} \, ,\qquad \quad
 - \widetilde{\sigma}_{11,2} = \widetilde{\sigma}_{33,2} \, , \qquad \quad
 \widetilde{\sigma}_{22,1} = - \widetilde{\sigma}_{33,1} \, .
\end{align}
Each condition in eq.\eqref{sym_m_from_Sigma} states that symmetric couple stress $\f m$
originates from just the spatial variation of deviatoric stress. It seems to be a physically reasonable constitutive requirement to distinguish couple stress from the gradient of dilational stress $\fg \sigma^d$. If we want couple stress to be independent of spatial pressure variations, the couple stress tensor must be chosen symmetric! Vice versa, if we want pressure gradients independent of couple stress we must choose the couple stress tensor $\f m$ to be symmetric.\footnote{Interestingly, Hadjesfandiari and Dargush \cite{hadjesfandiari2014evo}, coming from couple stress models for fluid dynamics, connect pressure gradients to couple stresses and assume therefore the opposite, namely that $\f m \in \so(3)$ is skew-symmetric.}

\begin{figure}
\centering
\includegraphics[height=30mm]{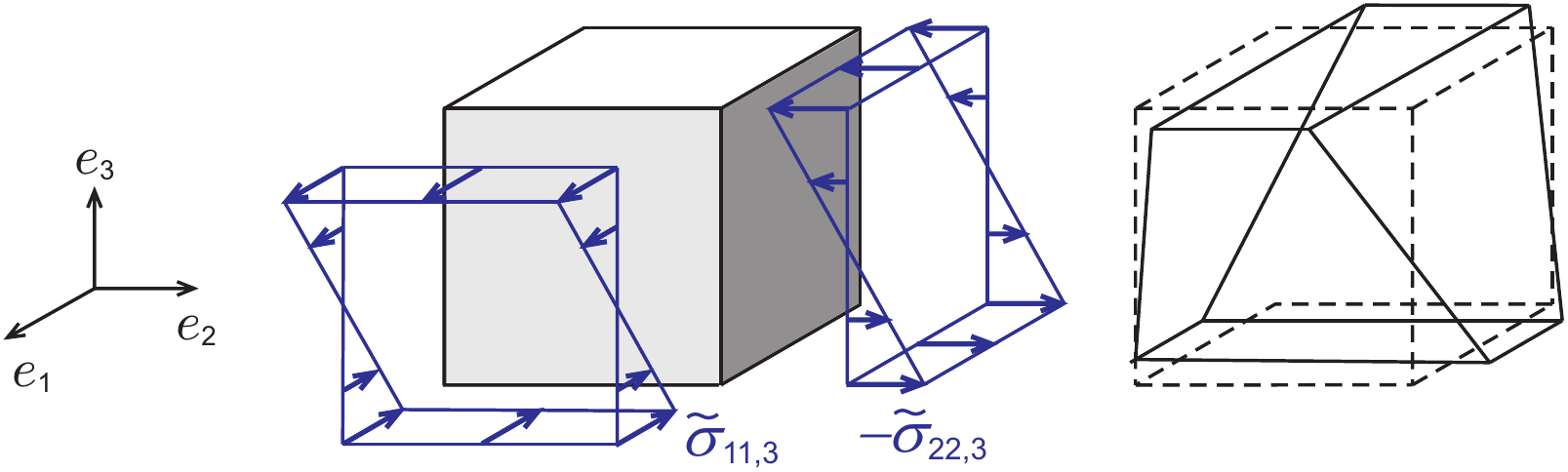}
\caption{Linear gradients of normal stress stretch and squeeze a finite cube from the symmetry condition of $\f m$ such that its overall volume is not affected.}
 \label{StretchSqueeze}
\end{figure}

Although the condition in eq.\eqref{sym_m_from_Sigma} concerns normal components of the stress tensor, it forces the stress field to be altered such that a deviatoric process results, as exemplarily drawn in
Fig.\ref{StretchSqueeze} for the condition given in eq.\eqref{sym_m_from_Sigma}$_1$. The
infinitesimal cube becomes stretched and squeezed such that alternating pure shear in the $e_1$-$e_2$- plane appears. Since we discuss here the second order stress resultant $\f m$, it is natural that the three conditions in eq.(\ref{sym_m_from_Sigma}) apply to \textit{gradients} of the first order quantity, namely the total stress $\widetilde{\fg \sigma}$.

\begin{remark} The assumption of symmetric couple stress tensors $\f m$ excludes that the change of volume is connected to couple stress. This is a physically meaningful assumption to decouple volumetric and deviatoric effects of secondary order within an elastic continuum theory of higher order.
\end{remark}

Some stress gradients in $\TDp.\Deltax$ \emph{do not} contribute to the couple stress tensor $\f m$. Such stress gradients appear from tangential tractions (shear components) varying parallel to their direction of action. This motivates to enhance our split from eq.\eqref{SplitSigma} further into
\begin{align}\label{SplitSigma2}
\TD = \TDnp + \underbrace{\TDpone + \TDptwo}_{\TDp} \, ,
\end{align}
with
\begin{align}\label{LinStressPolar2}
\TDpone . \Deltax =\!\!
   \left( \begin{array}{ccc}
     0 & \widetilde{\sigma}_{12,1} & \widetilde{\sigma}_{13,1} \\
     0 & 0 & 0 \\
     0 & 0 & 0 \\
   \end{array}
 \right) \Delta x_1
   +\!\left( \begin{array}{ccc}
    0 & 0 & 0 \\
    \widetilde{\sigma}_{21,2} & 0 & \widetilde{\sigma}_{23,2} \\
    0 & 0 & 0 \\
   \end{array}
 \right) \Delta x_2
   +\!\left( \begin{array}{ccc}
   0 & 0 & 0 \\
   0 & 0 & 0 \\
   \widetilde{\sigma}_{31,3} & \widetilde{\sigma}_{32,3} & 0 \\
   \end{array}
 \right) \Delta x_3 \,.
\end{align}
Since terms from $\TDpone . \Deltax$ do not correspond to couple stress, they may be related to conformal mappings, which are generally not isochoric. \\
\\
Next, we investigate the term $\TDb . \Deltaxb$. Similar to eq.\eqref{SplitSigma2} we
split
\begin{align}
\TDb = \TDbone + \TDbtwo \, ,
\end{align}
with
\begin{align}\label{BiLinStressNotPolar}
\TDbone . \Deltaxb =
   \left( \begin{array}{ccc}
     0 & \widetilde{\sigma}_{12,12} & \widetilde{\sigma}_{13,12} \\
     \widetilde{\sigma}_{21,12} & 0 & \widetilde{\sigma}_{23,12} \\
     0 & 0 & \widetilde{\sigma}_{33,12} \\
   \end{array}
 \right) & \Delta x_1 \, \Delta x_2
   +\left( \begin{array}{ccc}
    0 & \widetilde{\sigma}_{12,13} & \widetilde{\sigma}_{13,13} \\
    0 & \widetilde{\sigma}_{22,13} & 0 \\
    \widetilde{\sigma}_{31,13} & \widetilde{\sigma}_{32,13} & 0 \\
   \end{array}
 \right) \Delta x_1 \, \Delta x_3 \notag \\
   &+\left( \begin{array}{ccc}
   \widetilde{\sigma}_{11,23} & 0 & 0 \\
   \widetilde{\sigma}_{21,23} & 0 & \widetilde{\sigma}_{23,23} \\
   \widetilde{\sigma}_{31,23} & \widetilde{\sigma}_{32,23} & 0 \\
   \end{array}
 \right) \Delta x_2 \, \Delta x_3 \,,
\end{align}
and
\begin{align}\label{BiLinStressPolar}
\TDbtwo . \Deltaxb =
   \left( \begin{array}{ccc}
     \widetilde{\sigma}_{11,12} & 0 & 0 \\
     0 & \widetilde{\sigma}_{22,12} & 0 \\
     \widetilde{\sigma}_{31,12} & \widetilde{\sigma}_{32,12} & 0 \\
   \end{array}
 \right) & \Delta x_1 \, \Delta x_2
   +\left( \begin{array}{ccc}
    \widetilde{\sigma}_{11,13} & 0 & 0 \\
    \widetilde{\sigma}_{21,13} & 0 & \widetilde{\sigma}_{23,13} \\
    0 & 0 & \widetilde{\sigma}_{33,13} \\
   \end{array}
 \right) \Delta x_1 \, \Delta x_3 \notag \\
   &+\left( \begin{array}{ccc}
   0 & \widetilde{\sigma}_{12,23} & \widetilde{\sigma}_{13,23} \\
   0 & \widetilde{\sigma}_{22,23} & 0 \\
   0 & 0 & \widetilde{\sigma}_{33,23} \\
   \end{array}
 \right) \Delta x_2 \, \Delta x_3 \,.
\end{align}
This split is motivated by the equalities
\begin{align}\label{Lin_2a_Nonpolar}
\int_{\partial \B^i_c} \f r_i \times (\TDbone . \Deltaxb)  \cdot \f n_i \di A = 0 \,,
 \quad \quad i=1,...,6 \,
\end{align}
and
\begin{align}\label{Lin_2a_Nonpolarb}
\sum_{i=1}^6 \int_{\partial \B^i_c} \f x_P \times (\TDbone . \Deltaxb) \cdot \f n_i \di A =
0 \,,
\end{align}
stating that $\TDbone . \Deltaxb$ is nonpolar. However, we find
\begin{align}\label{BiLin_p_polar}
\int_{\partial \B^1_c}  \f r_1 \times (\TDbtwo . \Deltaxb)  \cdot \f n_1 \di A = \frac{L_c^5}{24}\left(
  \begin{array}{c}
    \widetilde{\sigma}_{31,21}- \widetilde{\sigma}_{21,31} \\
    \widetilde{\sigma}_{11,31} \\
    - \widetilde{\sigma}_{11,21} \\
  \end{array}
\right) \, , \notag \\
\int_{\partial \B^2_c}  \f r_2 \times (\TDbtwo . \Deltaxb)  \cdot \f n_2 \di A = \frac{L_c^5}{24}\left(
  \begin{array}{c}
    - \widetilde{\sigma}_{22,32} \\
    \widetilde{\sigma}_{12,32}- \widetilde{\sigma}_{32,12} \\
    \widetilde{\sigma}_{22,12} \\
  \end{array}
\right) \, , \notag \\
\int_{\partial \B^3_c}  \f r_3 \times (\TDbtwo . \Deltaxb)  \cdot \f n_3 \di A = \frac{L_c^5}{24}\left(
  \begin{array}{c}
    \widetilde{\sigma}_{33,23} \\
    - \widetilde{\sigma}_{33,13} \\
    \widetilde{\sigma}_{23,13}- \widetilde{\sigma}_{13,23} \\
  \end{array}
\right) \, , \\
\int_{\partial \B^4_c}  \f r_4 \times (\TDbtwo . \Deltaxb)  \cdot \f n_4 \di A =
 \int_{\partial \B^1_c}  \f r_1 \times (\TDbtwo . \Deltaxb)  \cdot \f n_1 \di A \, , \notag \\
\int_{\partial \B^5_c}  \f r_5 \times (\TDbtwo . \Deltaxb)  \cdot \f n_5 \di A =
\int_{\partial \B^2_c}  \f r_2 \times (\TDbtwo . \Deltaxb)  \cdot \f n_2 \di A \, , \notag \\
\int_{\partial \B^6_c} \f r_6 \times (\TDbtwo . \Deltaxb)  \cdot \f n_6 \di A =
\int_{\partial \B^3_c}  \f r_3 \times (\TDbtwo . \Deltaxb)  \cdot \f n_3 \di A \,,\notag
\end{align}
and
\begin{align}\label{BiLin_Bipolar}
\sum_{i=1}^6 \int_{\partial \B^i_c}  \f x_P \times (\TDbtwo . \Deltaxb)  \cdot \f n_i \di A
& = \frac{L_c^5}{12}\left(
  \begin{array}{c}
    \widetilde{\sigma}_{31,21}- \widetilde{\sigma}_{21,31}- \widetilde{\sigma}_{22,32}+ \widetilde{\sigma}_{33,23} \\
    \widetilde{\sigma}_{11,31}+ \widetilde{\sigma}_{12,32}- \widetilde{\sigma}_{32,12}- \widetilde{\sigma}_{33,13} \\
    - \widetilde{\sigma}_{11,21}+ \widetilde{\sigma}_{22,12}+ \widetilde{\sigma}_{23,13}- \widetilde{\sigma}_{13,23} \\
  \end{array}
\right) \notag \\
&= L_c^3\left(
  \begin{array}{c}
    m_{11,1}+m_{12,2}+m_{13,3} \\
    m_{21,1}+m_{22,2}+m_{23,3} \\
    m_{31,1}+m_{32,2}+m_{33,3} \\
  \end{array}
\right)= V_c \, \, \Div (\f m) \,.
\end{align}
\begin{remark}
Stress gradients $\TDbtwo . \Deltaxb$ are bipolar and are connected to the divergence of couple stress.\\
\end{remark}

\noindent
Finally, we investigate $\TDq . \Deltaxq$, obtaining
\begin{align}\label{Quad_notpolar}
\int_{\partial \B^i_c}  \f r_i \times (\TDq . \Deltaxq) \cdot \f n_i \di A = 0 \, , \qquad
i=1,...,6 \,,
\end{align}
and
\begin{align}\label{Quad_Semipolar}
\sum_i & \int_{\partial \B^i_c}  \f x_P \times  (\TDq . \Deltaxq)  \cdot \f n_i \di A \notag
\\ & = \frac{L_c^5}{12}\underbrace{\left(
  \begin{array}{c}
    \widetilde{\sigma}_{32,22}- \widetilde{\sigma}_{23,33} \\
    \widetilde{\sigma}_{13,33}- \widetilde{\sigma}_{31,11} \\
    \widetilde{\sigma}_{21,11}- \widetilde{\sigma}_{12,22} \\
  \end{array}
\right)}_{2 \, \axl \skew [\Grad{\widetilde{\fg \sigma}} : \overline{\fg \nabla}]} + \frac{L_c^5}{24} \underbrace{\left(
  \begin{array}{c}
    ( \widetilde{\sigma}_{32}- \widetilde{\sigma}_{23})_{,11}+( \widetilde{\sigma}_{32}- \widetilde{\sigma}_{23})_{,22}+( \widetilde{\sigma}_{32}- \widetilde{\sigma}_{23})_{,33}\\
    ( \widetilde{\sigma}_{13}- \widetilde{\sigma}_{31})_{,11}+( \widetilde{\sigma}_{13}- \widetilde{\sigma}_{31})_{,22}+( \widetilde{\sigma}_{13}- \widetilde{\sigma}_{31})_{,33}\\
    ( \widetilde{\sigma}_{21}- \widetilde{\sigma}_{12})_{,11}+( \widetilde{\sigma}_{21}- \widetilde{\sigma}_{12})_{,22}+( \widetilde{\sigma}_{21}- \widetilde{\sigma}_{12})_{,33}\\
  \end{array}
\right)}_{\Div \Grad{2 \, \axl \skew(\widetilde{\fg \sigma})}} \, .
\end{align}
Thus, $\TDq . \Deltaxq$ is semipolar.\\

\noindent
For compact symbolic notation in eq.\eqref{Quad_Semipolar} we define the nabla-operator $\overline{\fg \nabla} \in \R^{3\times 3\times 3}$ reading:
\begin{align}\label{NablaDreistufig}
\overline{\fg \nabla} \colonequals \nabla_k \, \f e_k \otimes \f e_k \otimes \f e_k \,.
\end{align}
The term $\Grad{\widetilde{\fg \sigma}} : \overline{\fg \nabla}$ in eq.\eqref{Quad_Semipolar} is extraordinary. It is similar to the divergence but of second order reading in index notation
\begin{align}\label{DoubleContraction3}
\Grad{\widetilde{\fg \sigma}} : \overline{\fg \nabla} &= \sigma_{ij,k} \, \f e_i \otimes \f e_j \otimes \f e_k : \nabla_k \, \f e_k \otimes \f e_k \otimes \f e_k = \sigma_{ij,k} \, \underbrace{\Scal{\f e_j \,,\, \f e_k}}_{\displaystyle = \delta_{jk}} \, \underbrace{\Scal{\f e_k \,,\, \f e_k}}_{\displaystyle =1} \, \nabla_k \, \f e_i \otimes \f e_k \notag \\ &=
\sigma_{ik,k} \, \nabla_k \, \f e_i \otimes \f e_k = \sigma_{ik,kk} \, \f e_i \otimes \f e_k \quad \in \R^{3\times 3} \,.
\end{align}
The first term in eq.\eqref{Quad_Semipolar} generally contributes to the balance of
angular momentum with
\begin{align}\label{Term1QuadMom}
 \frac{L_c^5}{12} \, 2 \, \axl \skew [\Grad{\widetilde{\fg \sigma}} : \overline{\fg \nabla}]
 = 2 \,V_c \,   \axl \skew \underbrace{\left[ \frac{L_c^2}{12} \Grad{\widetilde{\fg \sigma}} : \overline{\fg \nabla} \right]}_{\equalscolon \fg \chi}= 2 \,V_c \,  \axl \skew  [\fg \chi] \,.
\end{align}
The second term in eq.\eqref{Quad_Semipolar} contributes to the balance of angular
momentum if $\widetilde{\fg \sigma} \notin \Sym$, yielding
\begin{align}\label{Term2QuadMom}
 \frac{L_c^5}{24} \Div \Grad{2 \, \axl \skew(\widetilde{\fg \sigma})}
 = V_c \,  \Div \underbrace{\left[ \frac{L_c^2}{24} \, \Grad{2 \, \axl \skew(\widetilde{\fg \sigma})}
 \right]}_{\equalscolon \fg \psi} = V_c \, \, \Div [\fg \psi] \,.
\end{align}
\subsection{Balance of angular momentum}
Assuming a constant external loading from net couples $\f c$ within the cube $\B_c$, we
conclude from eq.\eqref{Sigma0Semipolar}, \eqref{BiLin_Bipolar}, and \eqref{Quad_Semipolar} that the balance of angular momentum reads
\begin{align}\label{Balanceofangularmomentum}
\sum_{i=1}^6 \int_{\partial \B^i_c}  \f x_P \times ( \widetilde{\fg \sigma}^0 + \TDbtwo . \Deltaxb + \TDq . \Deltaxq). \f n_i \di A + \int_{\B_c} \f c \di V & = 0 \notag \\
\Leftrightarrow V_c \, [ \Div \f m + \Div \fg \psi  + 2 \, \axl \skew ( \widetilde{\fg \sigma}^0 + \fg \chi ) + \f c ] &=0 \,.
\end{align}
We omit the mixture of polar and semipolar quantities in the Div-operator of
eq.\eqref{Balanceofangularmomentum}, because $\f m$ arises from a bipolar term of the
Taylor series expansion, whereas $\fg \psi$ originates from a semipolar term, which vanishes for $\widetilde{\fg \sigma} \in \Sym(3)$. However, couple stress is often introduced axiomatically via kinematic and constitutive assumptions such that both contributions to the balance of angular momentum may appear from a single quantity $\widetilde{\f m}= \f m + \fg \psi$. For $\widetilde{\fg \sigma} \in \Sym(3)$ the balance of angular momentum reduces to
\begin{align}\label{ReducedBoAM}
\Div \f m  + \fg \epsilon : \fg \chi ^T + \f c =0 \, , \qquad \tr (\f m) = 0 \,, \qquad \fg \chi = \frac{L_c^2}{12} \, \Grad{\widetilde{\fg \sigma}} : \overline{\fg \nabla}  \,.
\end{align}
In the appendix \ref{KapXsiObjectiv} we show, that $\fg \chi$ is objective, such that eq.\eqref{ReducedBoAM} is independent on the chosen basis.
\begin{remark}
Eq.\eqref{ReducedBoAM} is an extension to the balance of angular momentum in eq.\eqref{CoupleStressBalanceMomentum}. Thus, the linear couple stress theory (in this
interpretation) neglects $\TDq . \Deltaxq$. It is similar to the classical continuum
theory neglecting $\TDp . \Deltax$ in the balance of linear momentum.
\end{remark}

The couple stress vector $\f m_n = \f m \cdot \f n$ is {\bf a 1st moment}. Thus, it is independent of its point of application {\bf within a rigid body.} On first sight, assuming a rigid body to discuss the properties of couple stress seems to be allowed. Yang et al.~\cite{Yang02} mention that argument on page 2733: ``In conventional mechanics, a couple of forces is a free vector in the space of the material particle system. The couple can be translated and applied to any point in the system, which means that the motive effect of a couple on the system of material particles is independent of the location where the couple is applied. Thus, the forces $\f F_i$ and the couples of forces $\f L_i$ applied to a set of material particles within the system is equivalent to a resultant force and a resultant couple of forces, and the couple can be applied to an arbitrary point within the system." But rigidity and independence of point of application is irrelevant and even perplexing when discussing the properties of couple stress.\\

{\bf Force stress and couple stress localize to the center of infinitesimal cubes $\B_c$ to define balance equations.} Both force traction and couple traction arise from the same total stress function $\widetilde{\fg \sigma}$, its Taylor series expansion, and barycentric balance equations. Moving the point of application for couple stress would also move the stress function itself, which is not admissible. Thus, the properties of couple stress do not arise from an argument, which is only true in a rigid body. The properties of couple stress are polar properties of stress and its barycentric fluctuation.\\

Generally, deformation enters the continuum theory independent of balance equations via kinematic and constitutive equations. The derivation of static balance equations usually does not break down force stress by Taylor series expansion but also considers subdomains $\B_s$. By assuming constant net forces $\f f$, and constant net couples $\f c$ in $\B_s$, the balance of linear and angular momentum becomes
\begin{align}\label{MomentumBalance}
  \int_{\B_s}   \f f \, \di V + \int_{\partial \B_s}   \widetilde{\fg \sigma}_n(\f x) \, \di A =0,   \qquad \forall \, \B_s \subset \B \,
\end{align}
and
\begin{align}\label{AngularBalance}
  \int_{\B_s}  \f x \times   \f f +   \f c \, \di V + \int_{\partial \B_s}  \f x \times   \widetilde{\fg \sigma}_n(\f x) +   \f m_n(\f x) \, \di A =0 \qquad \forall \, \B_s \subset   \B \, ,
\end{align}
respectively. In eq.\eqref{AngularBalance} the couple traction $\f m_n(\f x)$ is axiomatic again and the position vector $\f x$ defines the lever arm of forces.
With the help of Cauchy's principle and the divergence theorem
\begin{align}\label{CauchyAndDiv1}
 \int_{\partial \B_s}  \widetilde{\fg \sigma}_n \, \di A =
 \int_{\partial \B_s}   \widetilde{\fg \sigma}  \cdot \f n \, \di A =
 \int_{\B_s}   \Div \widetilde{\fg \sigma} \, \di V \, ,
\end{align}
eq.\eqref{MomentumBalance} becomes
\begin{align}\label{MomBalVol}
\int_{\B_s}  ( \Div \widetilde{\fg \sigma}  + \f  f ) \, \di V =0 \qquad \forall \, \B_s
\subset \B \quad \Leftrightarrow \quad \Div \widetilde{\fg \sigma}  + \f  f  =0 \,,
\end{align}
which is in accordance with the local statement of eq.\eqref{Impulsbilanz1}$_3$. Similarly, using Cauchy's principle and the divergence theorem for the couple stress vector yields
\begin{align}\label{CauchyAndDiv2}
 \int_{\partial \B_s}   \f m_n \, \di A =
 \int_{\partial \B_s}   \f m  \cdot \f n \, \di A =
 \int_{\B_s}   \Div \f m \, \di V \, ,
\end{align}
as well as
\begin{align}\label{CauchyAndDiv3}\notag
 \int_{\partial \B_s}  &  \f x \times \widetilde{\fg \sigma}_n \, \di A=
 \int_{\partial \B_s}   \f x \times (\widetilde{\fg \sigma}.\f n) \, \di A \\
 &=\int_{\partial \B_s}   (\f x \times \widetilde{\fg \sigma}).\f n \, \di A =
 \int_{\B_s}   \Div (\f x \times \widetilde{\fg \sigma}) \, \di V =
 \int_{\B_s}  \f x \times \Div \widetilde{\fg \sigma} + 2\, \axl \skew\,\widetilde{\fg \sigma} \, \di V \, .
\end{align}
With eq.\eqref{CauchyAndDiv3} the balance of angular momentum in
eq.\eqref{AngularBalance} becomes
\begin{align}\label{AngBalVol}
  \int_{\B_s} \f x \times \left[ \underbrace{(\Div \widetilde{\fg \sigma} + \f f )}_{=0 \,   {\rm due \, to \, eq.\eqref{MomBalVol}}} + 2 \, \axl \skew\, \widetilde{\fg \sigma} +  \Div \f m + \f c \right] \, \di V & =0 \qquad \forall \, \B_s \subset \B \notag \\
\Leftrightarrow 2 \, \axl \skew\,\widetilde{\fg \sigma} + \Div \f m + \f c & =0 \,.
\end{align}
Putting together eq.\eqref{MomBalVol} and eq.\eqref{AngBalVol} we have obtained the
statement of balance equations:
\begin{align}\label{BalanceAgain}
\fbox{
\parbox[][0.2cm][c]{8cm}{$
\Div \widetilde{\fg \sigma}+ \f  f  = 0 \, , \qquad \Div \f m + 2 \, \axl
\skew\,\widetilde{\fg \sigma} + \f c  =0 \, , $}}
\end{align}
which are the force and moment balance laws governing the translational and rotational equilibrium by considering infinitesimal elements of matter and fully equivalent to system \eqref{CoupleStressBalanceMomentum}.\\

In order to augment the equations in the above box basing ourselves on particular constitutive relations we place ourselves in a small strain, isotropic linearized setting. There, the basic kinematical variables are the displacement gradient and we may {\bf constitutively} prescribe (only) the symmetric part of the total stress tensor by
\begin{align}\label{symsigma}
  \sym \widetilde{\fg \sigma} = 2 \, \mu \sym\text{Grad}[\f u] + \lambda \tr(\sym \text{Grad}[\f u]) \id \,.
\end{align}
Moreover,  consistent with isotropy we require the couple stress tensor $\f m$ to be an isotropic tensor function of the curvature tensor $\ks =\12 \, \Grad{\curl \f u} $. Then, the most general representation of $\f m$ is given by
\begin{align}\label{defmbycurl}
  \f m = \alpha_1 \, \dev \sym \ks  + \alpha_2 \, \skew \ks +
  \alpha_3  \underbrace{\tr (\ks)}_{= 0} \id \,.
\end{align}
There cannot be an independent constitutive prescription for $\skew \widetilde{\fg
\sigma}=\widetilde{\fg \tau}$, instead we have the requirement $\widetilde{\fg \tau} =- \12
\anti \Div \f m$, whence the name ``indeterminate couple stress model." Once the linear constitutive requirements (\ref{symsigma}), (\ref{defmbycurl}) are introduced, the resulting equations of equilibrium lose objectivity, as is already clear in linear elasticity.
\subsection{A symmetric stress function implying trace free couple stresses}\label{KapTraceFree}
The classical indeterminate couple stress model leads to trace free couple stresses, as seen in eq.\eqref{DefCoupleStress} and eq.\eqref{defmbycurl}. The question arises: Can we find a symmetric total force stress function $\widetilde{\fg \sigma}$ which is in accordance with this statement? Naturally, the normal components of the couple stress tensor $\f m$ appear as key figures. Thus, claiming couple stress to be skew-symmetric \cite{hadjesfandiari2014evo} would not allow for the following discussion.

Let us consider a cube $\B_m$ with dimensions $L_c$ and symmetric total force stress
$\widetilde{\fg \sigma}$ due to infinitesimal deformation. We are able to find a symmetric, linear stress function
\begin{align}\label{SymStressFunction}
\widetilde{\fg \sigma} =\widetilde{\fg \sigma}_a + \widetilde{\fg \sigma}_b +
\widetilde{\fg \sigma}_c = a \, \f B_a + b \, \f B_b + c \, \f B_c \, , \quad a,b,c \in \mathbb{R}, \quad \f B_a, \f B_b, \f B_c \in \Sym(3)
\end{align}
with off-diagonal basis elements
\begin{align}\label{Ba}
\f B_a = \left(
           \begin{array}{ccc}
             0 & -z & y \\
             -z & 0 & 0 \\
             y & 0 & 0 \\
           \end{array}
         \right) \, , \quad
\f B_b = \left(
           \begin{array}{ccc}
             0 & z & 0 \\
             z & 0 & -x \\
             0 & -x & 0 \\
           \end{array}
         \right) \, , \quad
\f B_c = \left(
           \begin{array}{ccc}
             0 & 0 & -y \\
             0 & 0 & x \\
             -y & x & 0 \\
           \end{array}
         \right) \, .
\end{align}
The origin of the orthogonal $x,y,z$ coordinate system is barycentric in $\B_m$ and
aligned to the $\f e_1$, $\f e_2$, $\f e_3$ directions, as drawn in Fig.\ref{CubeKOS}. Since the stress function from eq.\eqref{SymStressFunction} fulfills
\begin{align}\label{DivBaBbBc}
\Div \f B_a = 0 \, , \quad \quad \Div \f B_b = 0 \, , \quad \quad \Div \f B_c = 0 \,,
\end{align}
it satisfies the static balance of linear momentum in the absence of body forces $\f f$.

To satisfy balance of angular momentum, a classical continuum theory requires the stress function from eq.\eqref{SymStressFunction} with $a,b,c \rightarrow 0$ for any
infinitesimal body with $x \rightarrow \di x$, $y \rightarrow \di y$ and $z \rightarrow \di z$. In a couple stress theory with $x,y,z \in [-L_c/2,L_c/2]$
the stress function $\widetilde{\fg \sigma}$ yields the normal components of the couple stress tensor $\f m$ defined on $\partial \B_c$ by
\begin{align}\label{componentM11}
L_c^2 \, \, m_{11} = \left( \int_{-\frac{L_c}{2}}^{\frac{L_c}{2}}
\int_{-\frac{L_c}{2}}^{\frac{L_c}{2}} \left(\begin{array}{c}
   0 \\
   y \\
   z
   \end{array}\right)
\times \widetilde{\fg \sigma} \cdot \f e_1 \di y \, \di z \right) \cdot \f e_1
 \, ,
\end{align}
%
\begin{align}\label{componentM22}
L_c^2 \, \, m_{22} = \left( \int_{-\frac{L_c}{2}}^{\frac{L_c}{2}}
\int_{-\frac{L_c}{2}}^{\frac{L_c}{2}} \left(\begin{array}{c}
   x \\
   0 \\
   z
   \end{array}\right)
\times \widetilde{\fg \sigma} \cdot \f e_2 \di x \, \di z \right) \cdot \f e_2 \, ,
\end{align}
\begin{align}\label{componentM33}
L_c^2 \, \, m_{33} = \left( \int_{-\frac{L_c}{2}}^{\frac{L_c}{2}}
\int_{-\frac{L_c}{2}}^{\frac{L_c}{2}} \left(\begin{array}{c}
   x \\
   y \\
   0
   \end{array}\right)
\times \widetilde{\fg \sigma} \cdot \f e_3 \di x \, \di y \right) \cdot \f e_3 \,.
\end{align}
Evaluating eq.\eqref{componentM11} - \eqref{componentM33} for $\widetilde{\fg \sigma}_a$
yields
\begin{align}\label{componentsM112233a}
m_{11}^a=\frac{a}{L_c^2} \left( \int_{-\frac{L_c}{2}}^{\frac{L_c}{2}}
\int_{-\frac{L_c}{2}}^{\frac{L_c}{2}} \left(\begin{array}{c}
   y^2 + z^2 \\
   0 \\
   0
   \end{array}\right) \di y \, \di z \right) \cdot \f e_1
 = \frac{a \, J_p}{L_c^2} \, ,
\end{align}
\begin{align}\label{componentsM112233b}
m_{22}^a =\frac{a}{L_c^2} \left( \int_{-\frac{L_c}{2}}^{\frac{L_c}{2}}
\int_{-\frac{L_c}{2}}^{\frac{L_c}{2}} \left(\begin{array}{c}
   0 \\
   -z^2 \\
   0
   \end{array}\right) \di x \, \di z \right) \cdot \f e_2
   = - \frac{a \, J_2}{L_c^2} \, ,
\end{align}
\begin{align}\label{componentsM112233c}
m_{33}^a=\frac{a}{L_c^2} \left( \int_{-\frac{L_c}{2}}^{\frac{L_c}{2}}
\int_{-\frac{L_c}{2}}^{\frac{L_c}{2}} \left(\begin{array}{c}
   0 \\
   0 \\
   -y^2
   \end{array}\right) \di x \, \di y \right) \cdot \f e_3 = - \frac{a \, J_3}{L_c^2} \,,
\end{align}
with the polar moment of inertia $J_p=L_c^4/6$ and the moment of inertia $J_2=L_c^4/12$ and $J_3=L_c^4/12$, respectively. Equivalently, for $\widetilde{\fg \sigma}_b$ and $\widetilde{\fg \sigma}_c$ we obtain
\begin{align}\label{componentsM112233d}
m_{11}^b = - \frac{b \, J_1}{L_c^2} \, , \quad  \quad m_{22}^b = \frac{b \, J_p}{L_c^2} \, , \quad  \quad  m_{33}^b = - \frac{b \, J_3}{L_c^2} \,,
\end{align}
and
\begin{align}\label{componentsM112233e}
m_{11}^c = - \frac{c \, J_1}{L_c^2} \, , \quad  \quad m_{22}^c = - \frac{c \, J_2}{L_c^2} \, , \quad  \quad m_{33}^c = \frac{c \, J_p}{L_c^2} \,,
\end{align}
with $J_1=L_c^4/12$. In summary, the stress function $\widetilde{\fg \sigma}$ from
eq.\eqref{SymStressFunction} leads to the couple stress tensor
\begin{align}\label{CSTTraceFree}
\f m = \frac{L_c^2}{12} \left(
  \begin{array}{ccc}
    2a-b-c & 0 & 0 \\
    0 & 2b-a-c & 0 \\
    0 & 0 & 2c-a-b \\
  \end{array}
\right) \,,
\end{align}
which is in accordance with eq.\eqref{Def_m_from_Sigma}. Further, this couple stress tensor $\f m$ is symmetric and trace free for any choice of $a,b,c \in \mathbb{R}$.
\section{The gap in the initial motivation of the modified couple stress
theory}\label{KapYang}\setcounter{equation}{0}
Yang et al.~\cite{Yang02} define the residual body couple vector
\begin{align}
2 \, \axl \skew\,\widetilde{\fg \sigma} + \f c=:\f m^*
\end{align}
from body couples $\f c$ and skew-symmetric parts of the total stress tensor
$\widetilde{\fg \sigma}$. These quantities balance in case of a local continuum theory assuming $\f m^* \equiv 0$. Then, the total stress tensor $\widetilde{\fg \sigma}$ becomes symmetric if body couples $\f c$ are absent. This is the Cauchy-Boltzmann axiom:
\begin{align}\label{CauchyBoltzmann}
\fbox{
\parbox[][0.2cm][c]{6cm}{$
\widetilde{\fg \sigma} = \widetilde{\fg \sigma}^T \quad \text{Cauchy-Boltzmann Axiom} $}}
\end{align}
Yang et al. also define the cross product of the position vector $\f x$ with a couple $\f L$ as ``couple of couple" or ``moment of couple"\footnote{We denote a force couple by $\f L$ in accordance with the notation in Yang et al. The couple stress vector $\f m$ relates to an area $a$, such that $\f L = a \, \f m$.}
\begin{align}\label{YangCoupleMoment}
\f M \colonequals \f x \times \f L \, .
\end{align}
Thus, they presume individual points of application $\f x_i$ for couples $\f L_i$.
Neither the definition in eq.\eqref{YangCoupleMoment} nor the presumption of given points of application for couples is problematic. However, in a system of material particles the set of balance equations
\begin{align}\label{YangBalances}
\underbrace{\sum \f F_i = 0}_{\text{balance of linear momentum}} \,, \qquad
\underbrace{\sum(\f x_i \times \f F + \f L_i) = 0}_{\text{balance of angular momentum}} \, , \qquad {\color{red}{\underbrace{\sum \f x_i \times \f L_i=0}_{\text{proposed balance equation from Yang et al.}}}}
\end{align}
given in \cite{Yang02}, p.2736, is untenable concerning the third statement (marked in red)\footnote{The first two statements are clear and express the linear and angular
momentum.}. The statement in \cite{Yang02}, pp.2735 (above eq.(18)) initiates the
fallacy: ``The couple vector $\f L_A$ at $A$ in a system of material particles in Fig.2(a) is equivalent to a couple $\f L'_A$ and a couple of couples $\f M'_A$ applied to the point $B$ in Fig.2(c)." However, such an equivalence does not exist. It is motivated by the physical property of forces, while $\f L_A$ is a couple vector.

Some sentences above, in the same section, we can read: ``The couple of forces is a free vector in the conventional mechanics, which means that the effect of the couple applied on an arbitrary point in the space of the system of material particles is independent of the position of the point. In other words, the couple can translate to any point in space freely and the resulting motive effects are unchanged."

Yang et al. argue that a couple stress theory locates the point of application for
couples. We agree with this statement, but eq.\eqref{YangBalances}$_3$ is not a proper
balance equation even for rigid bodies\footnote{Lazopoulos \cite{Lazo09} and
Hadjesfandiari and Dargush \cite{hadjesfandiari2014evo} have also noted the
inappropriateness of Yang et al.'s additional balance equation.}. Thus, the statement in
eq.\eqref{YangBalances}$_3$ can yield no proof for the symmetry of couple stresses.\\

\begin{figure}
\centering
\includegraphics[height=45mm]{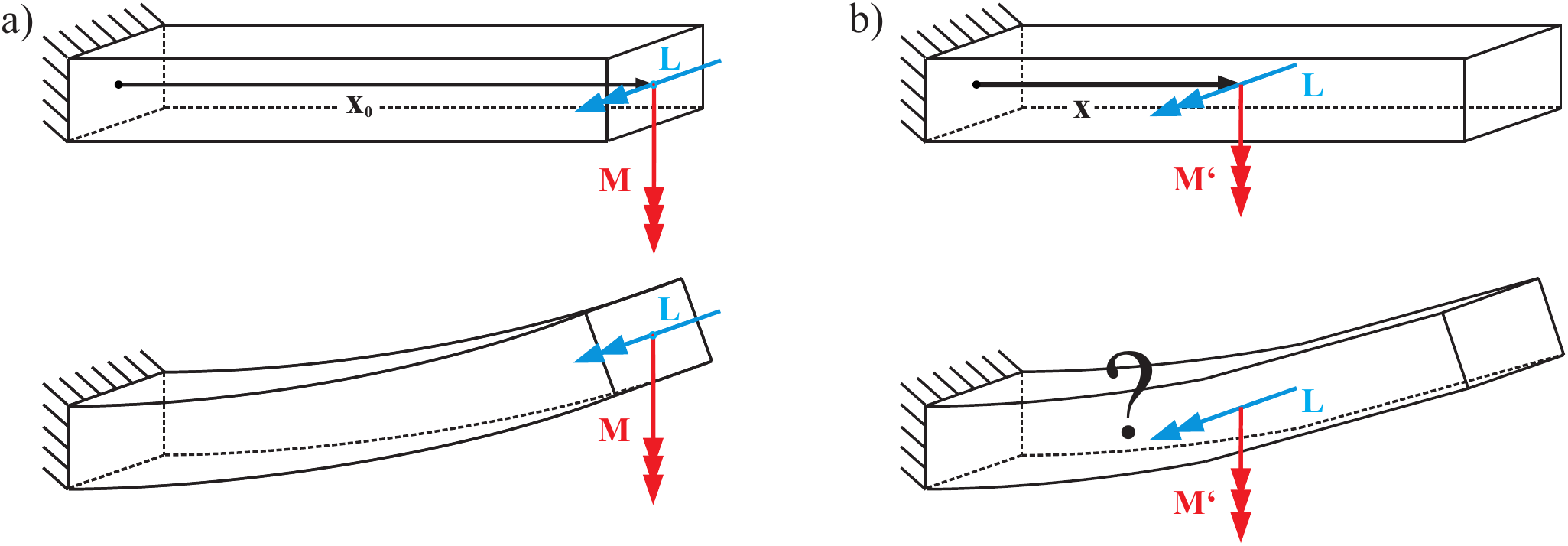}
\caption{a) Cantilever with couple $\f L$ applying at $\f x_0$. The couple generates bending within the whole structure and the moment of the couple reads $\f M = \f x_0 \times \f L$. b) Placing the couple to $\f x \neq \f x_0$  modifies the moment of the couple into $\f M' = \f x \times \f L$. Thus, the cantilever does not bend between $\f x$ and $\f x_0$. A compensating effect via $\f M'$ is questionable.}
 \label{Cantilever}
\end{figure}

The error occuring in eq.\eqref{YangBalances}$_3$ can be revealed in basic examples: consider e.g.~a cantilever with a couple vector $\f L$ as loading at its tip, see Fig.\ref{Cantilever}a. A basic choice for the origin of $\f x_0$ is the point of clamping. Then reaction forces and couple do not contribute to the sum in eq.\eqref{YangBalances}$_3$ since $\f x = 0$ at the point of clamping. But the ``moment of couple" $\f M=\f x_0 \times \f L$ does not vanish, since $\f x$ and $\f L$ are not linearly dependent. Hence, we find a simple example where eq.\eqref{YangBalances}$_3$ does not hold. This is independent on whether the cantilever is rigid or not.\\

Regarding in addition the cantilevers elasticity, the point of application for the couple vector $\f L$ is fundamental for its deformation. In Fig.\ref{Cantilever}b we placed the couple vector $\f L$ to the point $\f x \neq \f x_0$. From $\f x$ to the tip of the cantilever no curvature can appear, whereas in Fig.\ref{Cantilever}a it does. The ``moment of couple" $\f M=\f x \times \f L$ applies in $\f x$ but cannot compensate the lack of curvature between $\f x$ and the tip. Thus, the localization of couple vectors is not an exclusive requirement of couple stress theories but of any elastic theory. The shift of couple vectors from their point of application to another place is not allowed and can not be compensated by ``moments of couples". This is in accordance with the effect of forces in elasticity. It is not allowed to shift forces along their direction of action if the body is elastic\footnote{For the equilibrium of rigid bodies it is allowed to shift forces along their direction of action.}. Finally we remark that ``moments of couples" are not objective, since they depend on the origin of coordinates.

Nevertheless, Yang et al.~consider a generalization\footnote{Here, generalization means, that a statement motivated by point mechanics is used for a similar statement in continuum mechanics.} of eq.\eqref{YangBalances}$_3$ for all {\bf couple vectors} with their point of application $\f x$ and claim
\begin{align}\label{ThirdBalance}\notag
  \int_{V}  \f x \times \f m^*  \, \di V
  + \int_{\partial V}   \f x \times   \f m_n \, \di A &= 0 \quad \Leftrightarrow \\
  \int_{V}  \f x \times ( \f c + 2\,  \axl(\skew\, \widetilde{\fg \sigma}))  \, \di V
  + \int_{\partial V}  \f x \times   \f m_n \, \di A &=0 \qquad \text{for all subdomains} \qquad V\subset \B.
\end{align}
From Cauchy's principle and the divergence theorem for $\f m_n$ it follows that
\begin{align}\label{CauchyAndDiv4}
 \int_{\partial V}   \f x \times \f m_n \, \di A =
 \int_{\partial V}   (\f x \times \f m) \cdot \f n \, & \di A =
 \int_{V}   \Div (\f x \times \f m) \, \di V= \int_{V}  \f x \times \Div \f m + 2 \,  \axl(\skew\, \f m) \, \di V.
\end{align}
Thus, eq.\eqref{ThirdBalance} can be rewritten as
\begin{align}\label{ThirdBalanceVol}
  \int_{V}  \left[\f x \times (\underbrace{2\,  \axl(\skew\,\widetilde{\fg \sigma}) + \Div \f m + \f c}_{=0\ \text{from \eqref{AngBalVol}}} )
  + 2\,  \axl(\skew \, \f m)\right]\, \di V=0.
\end{align}
Since the balance of angular momentum is given by eq.\eqref{AngBalVol}, Yang et al.~conclude
\begin{align}
  \int_{V}  2\,  \axl(\skew \, \f m)\, \di V=0 \qquad \text{for all subdomains} \qquad V\subset \B.
\end{align}
Assuming that $\f m$ is a continuous couple stress tensor field it follows by localization that
\begin{align}\label{SymMomStress}
   \axl(\skew\, \f m) = 0 \quad \Leftrightarrow \quad \skew \, \f m =0
\end{align}
i.e.~that the couple stress tensor $\f m$ must be symmetric. However, we do not agree with eq.\eqref{ThirdBalance}, since the cross product of couple vectors with arbitrary position vectors does not vanish - except if position vectors and couple vectors are linearly dependent, see also the appendix \ref{KapPosIndependentM}. But such a linear dependency is arbitrary and not a physical law. Thus, the argument by Yang et al.~is incomplete, even though as seen in the previous sections, a symmetric couple stress tensor $\f m$ is indicated on different grounds.
\section{Torsion example}\label{KapTorsionExample}
\begin{figure}
\centering
\includegraphics[height=50mm]{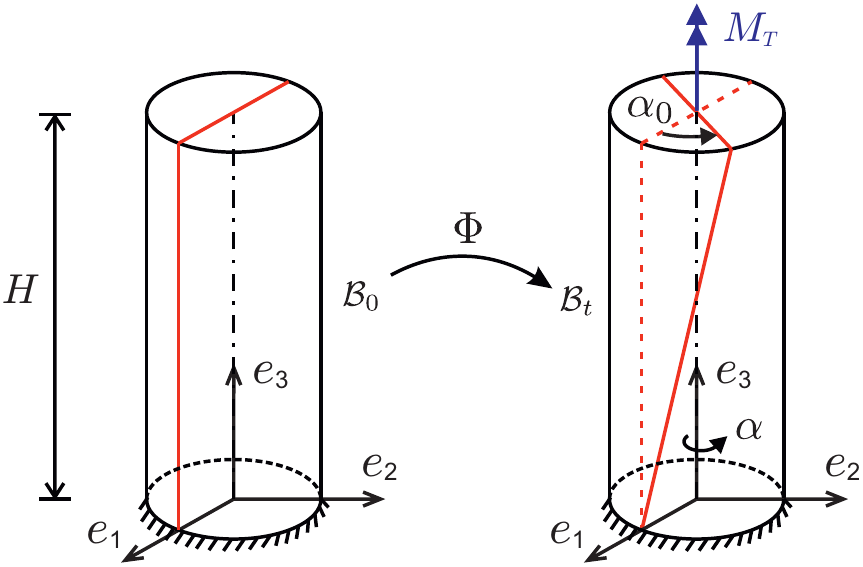}
\caption{Clamped circular beam obtaining simple torsional deformation $\fg \Phi$ from the moment $M_T$ at its tip.}
 \label{TorsionExampleZylinder}
\end{figure}

\begin{figure}
\centering
\includegraphics[height=50mm]{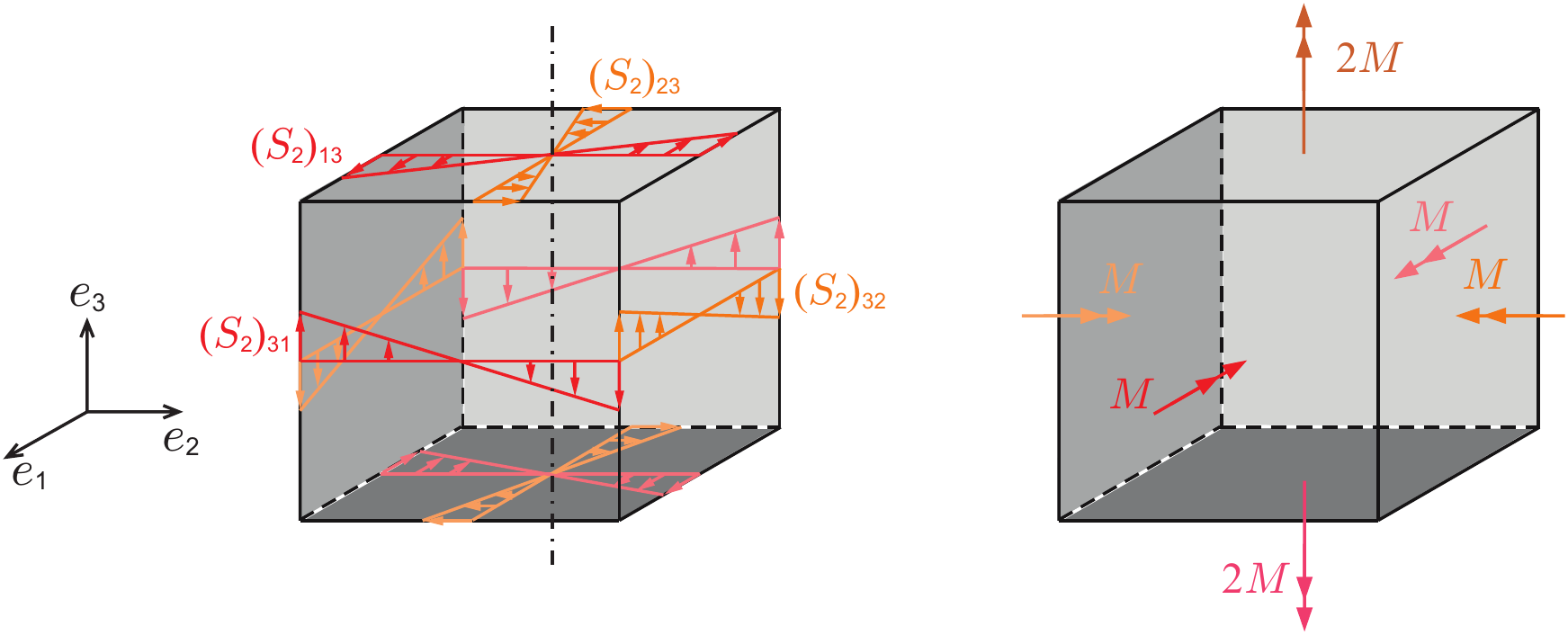}
\caption{a) Trend of components of the second Piola-Kirchhoff-stress tensor $\f S_2$ for an axially centered cube. b) Couples generated by tractions from  $\f S_2$ concerning the center of corresponding faces.}
 \label{CubeMom3Equilibri}
\end{figure}

\setcounter{equation}{0}
Let us consider simple torsional deformation of a circular beam by the angle $\alpha$ as linear function of the Lagrangian coordinate $z$, reading
\begin{align}
\alpha(z)=\frac{\alpha_0}{H}\,z = \overline{\alpha}\, z \, , \qquad \f x= \left(\begin{array}{c}x\\ y\\ z\\\end{array}\right) \, , \qquad \f x \in \B_0 \,.
\end{align}
The beam is fully clamped at $z=0$ and loaded by the moment $\f M_T$ in $\f e_3$-direction at its tip $z=H$. Each horizontal cross section rotates uniformly by the angle $\alpha$ with $\alpha_0$ at $z=H$. Thus, we define the constant gradient of rotation by $\overline{\alpha} = \alpha_0 / H$. In Fig.\ref{TorsionExampleZylinder} the system is shown in its initial state $\B_0$ and after deformation indicated by the actual state $\B_t$.
In the actual state the Lagrangian position vector $\f x$ is mapped to $\xbar$ via the rotation $\f R(z)$ reading
\begin{align}
\xbar = \f R \, \f x \, , \qquad \text{with} \quad \f R= \left[\begin{array}{ccc}
\cos(\overline{\alpha}\, z) & -\sin(\overline{\alpha}\, z) & 0\\
\sin(\overline{\alpha}\, z) & \cos(\overline{\alpha}\, z) & 0\\
0 & 0 & 1
\end{array}\right] \in \text{SO(3)} \,.
\end{align}
Abbreviating $\cos(\overline{\alpha}\, z)=c \, , \, \sin(\overline{\alpha}\, z)=s$, the actual position vector $\xbar$, the displacement vector $\f u$, and its gradient read
\begin{align}
& \xbar= \left(\begin{array}{c} c\,x-s\,y \\ s\,x+c\,y \\ z\end{array}\right) \, , \quad
\f u=\xbar-\f x=\left(\begin{array}{c} (c-1)\,x-s\,y \\ s\,x+(c-1)\,y \\ 0 \\\end{array}\right) \, , \notag \\
&\Grad{\f u}= \left[\begin{array}{ccc}
(c-1)& -s & (-s\,x-c\,y)\overline{\alpha}\\
s & (c-1) & (c\,x-s\,y)\overline{\alpha}\\
0 & 0 & 0\\
\end{array}\right]
 \,.
\end{align}
Simple torsion does not alter the position of points in the direction of the beam, which is the $\f e_3$-axis here. Further, it is an isochoric deformation, resulting from the determinant of the deformation gradient $\f F = \Lin \fg \Phi$ given by
\begin{align}
\f F= \id + \Grad{\f u}=\left[\begin{array}{ccc}
c& -s & (-s\,x-c\,y)\overline{\alpha}\\
s & c & (c\,x-s\,y)\overline{\alpha}\\
0 & 0 & 1\\
\end{array}\right] \, , \quad
\det \,\f F = 1 \cdot \det \left[\begin{array}{cc}
c & -s\\
s & c
\end{array}\right]
=c^2+s^2=1 \, .
\end{align}
 Next, let us assume isotropic and elastic Saint-Venant-Kirchhoff material. Therefore, we calculate the Green-Lagrange strain tensor:
\begin{align}\label{TorsionGreen}
\f E &= \12(\f F^T \, \f F - \id) \notag \\
&=\12 \, \left[\begin{array}{ccc}
c^2+s^2 -1 & 0 & -\overline{\alpha}(s\,c\,x+c^2\,y-s\,c\,x+s^2\,y)\\[2mm]
 & s^2+c^2 -1 & \overline{\alpha}(s^2\,x+s\,c\,y+c^2\,x-s\,c\,y)\\[2mm]
{\rm sym} &  &
\overline{\alpha}^2(s^2\,x^2+2\,s\,c\,x\,y+c^2\,y^2+c^2\,x^2-2\,s\,c\,x\,y+s^2\,y^2)
\end{array}\right] \notag \\
&= \12 \, \left[\begin{array}{ccc}
0 & 0 & -\overline{\alpha} \, y\\
 & 0 & \overline{\alpha} \, x\\
{\rm sym} &  & \underbrace{\overline{\alpha}^2 \, (x^2+y^2)}_{\approx 0}\\
\end{array}\right] \,.
\end{align}
Since we want to discuss linear couple stress models here, we restrict this example to small rotations with $\overline{\alpha}^2 \approx 0$ yielding $ \tr \f E \approx 0$. Thus, the volumetric inner energy in this example vanishes and the second Piola-Kirchhoff-stress tensor $\f S_2$ becomes
\begin{align}
W^{\rm SVK}=\mu \, \vert\vert \f E \vert\vert^2 + \frac{\lambda}{2}(\underbrace{\tr \f E}_{\approx 0})^2 \qquad \Rightarrow \f S_2=2\mu\f E + \lambda \tr \f E \id \approx  \mu \left[\begin{array}{ccc}
0 & 0 & -\overline{\alpha}\,y\\
 & 0 & \overline{\alpha}\,x\\
{\rm sym} &  & 0\\
\end{array}\right] \,.
\end{align}
Since the components of $\f S_2$ are defined in $\B_0$, we can draw them in the reference state, which is not rotated. In Fig.\ref{CubeMom3Equilibri}a we consider an axially centered cube with mid point at $x=0$, $y=0$, and arbitrary $z$.
If we consider the dimension of the cube to be $\di x$, then, the maximum total value of stress components is given by
\begin{align}
\bar{S} = | \12 \, \overline{\alpha} \, \mu \,  \di x | \,.
\end{align}
Next, let us calculate the couple along $\f e_1$ on $\partial \B^1_c$ from the traction
\begin{align}\label{TractionOn1}
\f S_2  \cdot \f n_1 = \left(
                 \begin{array}{c}
                   (S_2)_{11} \\
                   (S_2)_{21} \\
                   (S_2)_{31} \\
                 \end{array}
               \right)= \left(
                 \begin{array}{c}
                   0 \\
                   0 \\
                   - \mu \, \overline{\alpha} \, y  \\
                 \end{array}
               \right) \, ,
\end{align}
yielding
\begin{align}\label{CoupleOn1}
\f M&=\int_{\partial \B^1_c}  \f r_1 \times \f S_2  \cdot \f n_1 \di A_1
=\int_{\partial \B^1_c}  \left(
                 \begin{array}{c}
                   0 \\
                   y \\
                   z  \\
                 \end{array}
               \right)  \times  \left(
                 \begin{array}{c}
                   0 \\
                   0 \\
                   - \mu \, \overline{\alpha} \, y  \\
                 \end{array}
               \right) \di A_1 =  \int_{-\frac{\di x}{2}}^{\frac{\di x}{2}}  \int_{-\frac{\di x}{2}}^{\frac{\di x}{2}} -\mu\,\overline{\alpha}\,y^2\,\di y \, \di z \, \f e_1 \notag \\
&= - \mu \, \overline{\alpha} \, \di x \left[\frac{ y^3 }{3} \right]_{-\frac{\di x}{2}}^{\frac{\di x}{2}}  \, \f e_1 = - \mu \, \overline{\alpha} \, \frac{\di x^4}{12}  \, \f e_1 \, .
\end{align}
Similarly, we can calculate the couples resulting from $(S_2)_{13},\, (S_2)_{23}$ and $(S_2)_{32}$ on the other surfaces of the cube. We obtain that couples in $\f e_3$-direction are doubled and reverse to the result in eq.\eqref{CoupleOn1}. The results are drawn in Fig.\ref{CubeMom3Equilibri}b. Since the whole set of couples are constant in $\B$ and in a state of self-equilibrium, the balance of angular momentum remains classical, stating $\skew \fg \sigma = 0$. Moreover, the Saint-Venant-Kirchhoff material does not attain such a constant state of inner couples with curvature energy. However, a strain gradient theory accounts for such couples by introducing a curvature measure and additional constitutive laws. How must we constitute the linear indeterminate couple stress theory to be in accordance with the stress state from above?\\

Naturally, the linear strain measure $\sym \Grad{\f u}$ is equivalent to the Green strain measure $\f E$ in eq.\eqref{TorsionGreen} for small rotations $\overline{\alpha} \ll 1$ assuming $\overline{\alpha}^2 =0$, $c \rightarrow 1$ and $s \rightarrow \overline{\alpha}\,z$:
\begin{align}
& \f u^{\rm lin}=\left(\begin{array}{c} -\overline{\alpha}\,y\,z\\ \overline{\alpha}\,x\,z\\0\end{array}\right) \, , \quad \Grad{\f u^{\rm lin}}=\left[\begin{array}{ccc}
0 & -\overline{\alpha}\,z & -\overline{\alpha}\,y\\
\overline{\alpha}\,z & 0 & \overline{\alpha}\,x\\
0 & 0 & 0\\
\end{array}\right] \, , \notag \\
& \sym \Grad{\f u^{\rm lin}}=\12 \, \left[\begin{array}{ccc}
0 & 0 & -\overline{\alpha}\,y\\
0 & 0 & \overline{\alpha}\,x\\
-\overline{\alpha}\,y & \overline{\alpha}\,x & 0\\
\end{array}\right] \,.
\end{align}
Thus, the linear stress tensor becomes equivalent to the second Piola-Kirchhoff-tensor $\f S_2$ for the above assumptions. Further, the curvature $\ks$ is given by
\begin{align}
\ks &= \12 \Grad{\curl \f u^{\rm lin}} = \Grad{ \axl \skew \Grad{\f u^{\rm lin}}}
= \rm{Grad} \left\{ \axl \left[\begin{array}{ccc}
0 & -\overline{\alpha}\,z & -\12\overline{\alpha}\,y\\
\overline{\alpha}\,z & 0 & \12\overline{\alpha}\,x\\
\12\overline{\alpha}\,y & -\12\overline{\alpha}\,x & 0\\
\end{array}\right] \right\} \notag \\
&= \rm{Grad} \left\{ \overline{\alpha} \left(\begin{array}{c} -\12x \\ -\12 y \\ z \end{array}\right) \right\}
= \overline{\alpha}\left[\begin{array}{ccc}
-\12 & 0 & 0\\
0 & -\12 & 0\\
0 & 0 & 1
\end{array}\right] \,.
\end{align}
Thus, the linear indeterminate couple stress model from Section \ref{KapIndeterminate} generates the couple stress tensor
\begin{align}
\f m = 2 \, \mu \, L_c^2 \, ( \alpha_1 \,\sym \f k + \alpha_2 \, \skew \f k) =\overline{\alpha}  \, \mu \, L_c^2 \,  \alpha_1 \, \left[\begin{array}{ccc}
-1 & 0 & 0\\
0 & -1 & 0\\
0 & 0 & 2
\end{array}\right] \, .
\end{align}
Comparing $\f m$ with the couple generated by the traction from eq.\eqref{TractionOn1} we obtain
\begin{align}
\int_{\partial B_1^c} m_{11} \,  \di A \, \f e_1 = \f M \quad \Leftrightarrow \quad
- \overline{\alpha}  \, \mu \, L_c^2 \,  \alpha_1 \, \di x^2 \, \f e_1 =  - \mu \, \overline{\alpha} \, \frac{\di x^4}{12}  \, \f e_1 \quad \Leftrightarrow \quad L_c^2 \, \alpha_1 = \frac{\di x^2}{12} \,.
\end{align}
Obviously, the internal length scale $L_c$ corresponds to the dimension $\di x$ of the underlying cube to set up the couple stress itself from gradients of stress. Considering $ L_c^2 \, \alpha_1 > 0$ yields the indeterminate couple stress model to become stiffer than the Saint-Venant-Kirchhoff material. In the limit case $ L_c^2 \, \alpha_1 \rightarrow \infty$ the external couple $M_T$ is balanced by the constant component $m_{33}$ of the couple stress tensor $\f m$ such that $\overline{\alpha} \rightarrow 0$ and one observes unbounded stiffness of the torsion beam.\\

The symmetric stress function in Section \ref{KapTraceFree} includes this example with parameters $a=0$, $b=0$ and $c= 12 \, \overline{\alpha} \, \mu \, \alpha_1$. Therefore, the general case $a,b,c \in \R$ in Section \ref{KapTraceFree} represents the arbitrary mode of spatial torsion.
\section{Conclusions and outlook}\label{KapConclusion}
Couple stresses are usually postulated for higher gradient theories from constitutive
assumptions on curvature energy or from kinematical considerations. In this paper, an
independent approach is used to discuss the properties of couple stresses within
continuum theories. We assume that couple stress can be represented in terms of stress gradients at a finite Cauchy cube, which is the basic model for mechanical
equilibrium equations. This is in accordance with usual conceptions that couple stress appears
\begin{itemize}
  \item in the vicinity of mechanical singularities (i.e. nooks) because of large stress gradients.
  \item at small or thin samples with bending or torsion deformation. Since curvature is size dependent, it increases its value on small scales, where appropriate
        stress gradients need to be transformed into couples.
  \item for material with distinct inner structure on small scales accounting for curvature as mentioned above (foams, granular material).
\end{itemize}
At the infinitesimal level, we show that stress gradients may be divided into distinct parts contributing to the balance of linear momentum and into several terms contributing to the balance of angular momentum. Since couple stresses are postulated to balance angular momentum, we can identify these terms arising from the Taylor series expansion of total force stress $\widetilde{\fg \sigma}$. Next, by postulating symmetry of the total stress $\widetilde{\fg \sigma}$, we find that the couple stress tensor $\f m$ must be traceless. Further, assuming isochoric deformation via couple stress, we find an argument for its symmetry. This is consistent with a proposed variant of the linear isotropic indeterminate couple stress model with symmetric local force-stress, symmetric non-local force-stress, symmetric couple-stresses and complete traction boundary conditions published recently \cite{GhiNeMaMue15}.

However, we do not agree with the argument used in Yang et al.~\cite{Yang02} for the intrinsic symmetry of the couple stress tensor. In accordance with Hadjesfandiari and Dargush \cite{hadjesfandiari2014evo} we question their symmetry argument, which is a physically artificial postulate. On the other hand, we also challenge the statement from Hadjesfandiari and Dargush \cite{hadjesfandiari2013skew} that the couple stress tensor $\f m$ is purely skew-symmetric: our development clearly shows the contrary under suitable hypotheses. In our point of view, symmetry of the couple stress tensor $\f m$ is a physically consistent additional constitutive requirement.\\

In \cite{MuenchWoehler2016,MuenchWoehler2017} we have used the Taylor expansion of stress up to order three together with the approach from Section \ref{KapTaylor}. The terms of order three do not enter the balance of angular momentum but extend the balance of linear momentum into
\begin{align}\label{ExtendedLinearMomentum}
\Div \fg \sigma + \f f - \dot{\rho} \, \dot{\f x} - \rho \, \ddot{\f x} + \dfrac{1}{24} \, L_c^2 \,\, \Div \Grad{\Div \fg \sigma + \f f - \dot{\rho} \, \dot{\f x} - \rho \, \ddot{\f x}} = 0 \,.
\end{align}
As discussed in \cite{MuenchWoehler2017}, to use the above equation needs motivation from strong heterogeneities in elastic or inertial properties of the material or loading.\footnote{An example in \cite{MuenchWoehler2017} shows, that eq.\eqref{ExtendedLinearMomentum} allows for an approximated stress field of a periphractic structure within the material.} Similarly, the extended balance of angular momentum in eq.(3.51) relates our approach to the indeterminate couple stress and Cosserat theory, respectively. Proceeding the Taylor expansion up to order four yields additional terms in the balance of angular momentum over again and will be discussed in a subsequent paper.\\

\noindent
{\bf Acknowledgement}\\

\noindent
The work of I.D. Ghiba was supported by a grant of the Romanian National Authority for Scientific Research and Innovation, CNCS-UEFISCDI, project number PN-II-RU-TE-2014-4-1109.

\addcontentsline{toc}{section}{Appendix}
\appendix
\setcounter{section}{1} \setcounter{equation}{0}
\subsection{Objectivity of $\fg \chi$}\label{KapXsiObjectiv}
Let $\f Q$ be a constant rotation tensor with $\f Q^T \, \f Q=\f Q \, \f Q^T = \id$ and $\det\f Q = +1$, mapping the orthogonal referential system of Euclidean vectors $\f e_i$ to a rotated system $\f d_i$ by
\begin{align}\label{QmapsEuklid}
\f d_i = \f Q \, \f e_i \, , \quad \f Q = \f d_i \otimes \f e_i \, , \quad Q_{ia} = \Scal{\f d_i , \f e_a}  = \Scal{\f e_a , \f d_i} \, , \quad i=1,2,3 \, , \quad a=1,2,3 \,.
\end{align}
The components $Q_{ia}$ are defined by the commutative inner product. However, the rotation tensor $\f Q$ is generally not symmetric: $\f Q \neq \f Q^T$. Since the referential basis vectors $\f e_i$ are considered to be orthogonal and of unit length, the rotated basis vectors $\f d_i$ are orthogonal and of unit length as well:
\begin{align}
\Scal{\f d_i , \f d_j} = \Scal{Q_{ia} \, \f e_a , Q_{jb} \ \f e_b} = Q_{ia} \, Q_{jb} \, \delta_{ab}  = Q_{ia} \, Q_{ja}  = Q_{ia} \, Q^T_{aj} = \delta_{ij}\, , \quad \Scal{\f d_i , \f d_i} = 1 \, .
\end{align}
A second order tensor, for instance the stress tensor $\fg \sigma = \sigma_{ij} \, \f e_i \otimes \f e_j$, is objective, if it is independent of the referential system. Thus, the tensor components $\sigma^{\sharp}_{mn}$ with rotated basis $\f d_m \otimes \f d_n$ need to transform by
\begin{align}\label{ObjectivitySigma}
\sigma^{\sharp}_{mn} \, \f d_m \otimes \f d_n &= \sigma_{ij} \, \f e_i \otimes \f e_j \notag \\
\Leftrightarrow \quad \sigma^{\sharp}_{mn} \, \underbrace{\Scal{\f d_m \, , \, \f d_m}}_{\displaystyle =1} \, \underbrace{\Scal{\f d_n \, , \, \f d_n}}_{\displaystyle =1} &= \sigma_{ij} \, \, \underbrace{\Scal{\f d_m \, , \, \f e_i}}_{\displaystyle =Q_{mi}} \, \underbrace{\Scal{\f e_j \, , \, \f d_n}}_{\displaystyle =Q_{nj}} \notag \\
\Leftrightarrow \quad \sigma^{\sharp}_{mn} &= Q_{mi} \, \sigma_{ij} \, Q^T_{jn}\, .
\end{align}
If the components of the second order tensor $\fg \chi$ transform like the components of $\fg \sigma$ in eq.\eqref{ObjectivitySigma}, it is objective.\\

\noindent
Derivatives concerning the position $\f x = x_i \, \f e_i$ in the basis $\f e_i$ are defined by the nabla operator
\begin{align}
\fg \nabla_{\f x}(\ldots) = \frac{\partial(\ldots)}{\partial \f x} = \frac{\partial(\ldots)}{\partial x_i} \, \f e_i = (\nabla_{x})_i \, \f e_i \, .
\end{align}
To define the nabla operator in the rotated referential system $\f d_j$, we introduce the rotated position vector\\
\noindent ${\fg \xi = \xi_i \, \f d_i = \f Q \cdot \f x = Q_{ij} \, x_j \, \f e_i}$, yielding
\begin{align}
\fg \nabla_{\fg \xi}(\ldots) &= \frac{\partial(\ldots)}{\partial \fg \xi} = \frac{\partial(\ldots)}{\partial \xi_j} \, \f d_j = \underbrace{\frac{\partial(\ldots)}{\partial x_i}}_{\displaystyle (\nabla_x)_i} \, \frac{\partial x_i}{\partial \xi_j} \, \f d_j = (\nabla_{x})_i \, \frac{\partial ( Q^T_{ia} \, \xi_a) }{\partial \xi_j} \, \f d_j \notag \\
&= (\nabla_{x})_i \, Q^T_{ia} \, \delta_{aj} \, \f d_j =Q_{ji} \, (\nabla_{x})_i \, \f d_j \, .
\end{align}
Using eq.\eqref{DoubleContraction3}, the components of the tensor $\fg \chi$ from eq.\eqref{Term1QuadMom} are given in the rotated referential system by
\begin{align}\label{ChiObjective}
\frac{12}{L_c^2} \, \fg \chi = \Grad{\fg \sigma} : \overline{\fg \nabla}_{\fg \xi}
= \chi^{\sharp}_{mn} \, \f d_m \otimes \f d_n &= \sigma^{\sharp}_{mn,nn} \,\f d_m \otimes \f d_n = \sigma^{\sharp}_{mn}  \, (\nabla_{\xi})_n \, (\nabla_{\xi})_n \,\f d_m \otimes \f d_n \notag \\
&= Q_{mi} \, \sigma_{ij} \, \underbrace{Q^T_{jn} \, Q_{na}}_{\displaystyle \delta_{ja}} \, (\nabla_{x})_a \, Q_{na} \, (\nabla_{x})_a \, \f d_m \otimes \f d_n \notag \\
&= Q_{mi} \, \underbrace{\sigma_{ia} \, (\nabla_{x})_a \, (\nabla_{x})_a}_{\displaystyle \sigma_{ia,aa} = \chi_{ia}}  \, Q^T_{an} \, \f d_m \otimes \f d_n \notag \\
&= Q_{mi} \, \chi_{ia} \, Q^T_{an} \, \f d_m \otimes \f d_n \,.
\end{align}
Since the components of $\fg \chi$ transform by $\chi^{\sharp}_{mn} =Q_{mi} \, \chi_{ia} \, Q^T_{an}$ it is objective.
\subsection{The classical balance of linear momentum}\label{KapClassicBalanceLinMom}
\begin{figure}
\centering
\includegraphics[height=80mm]{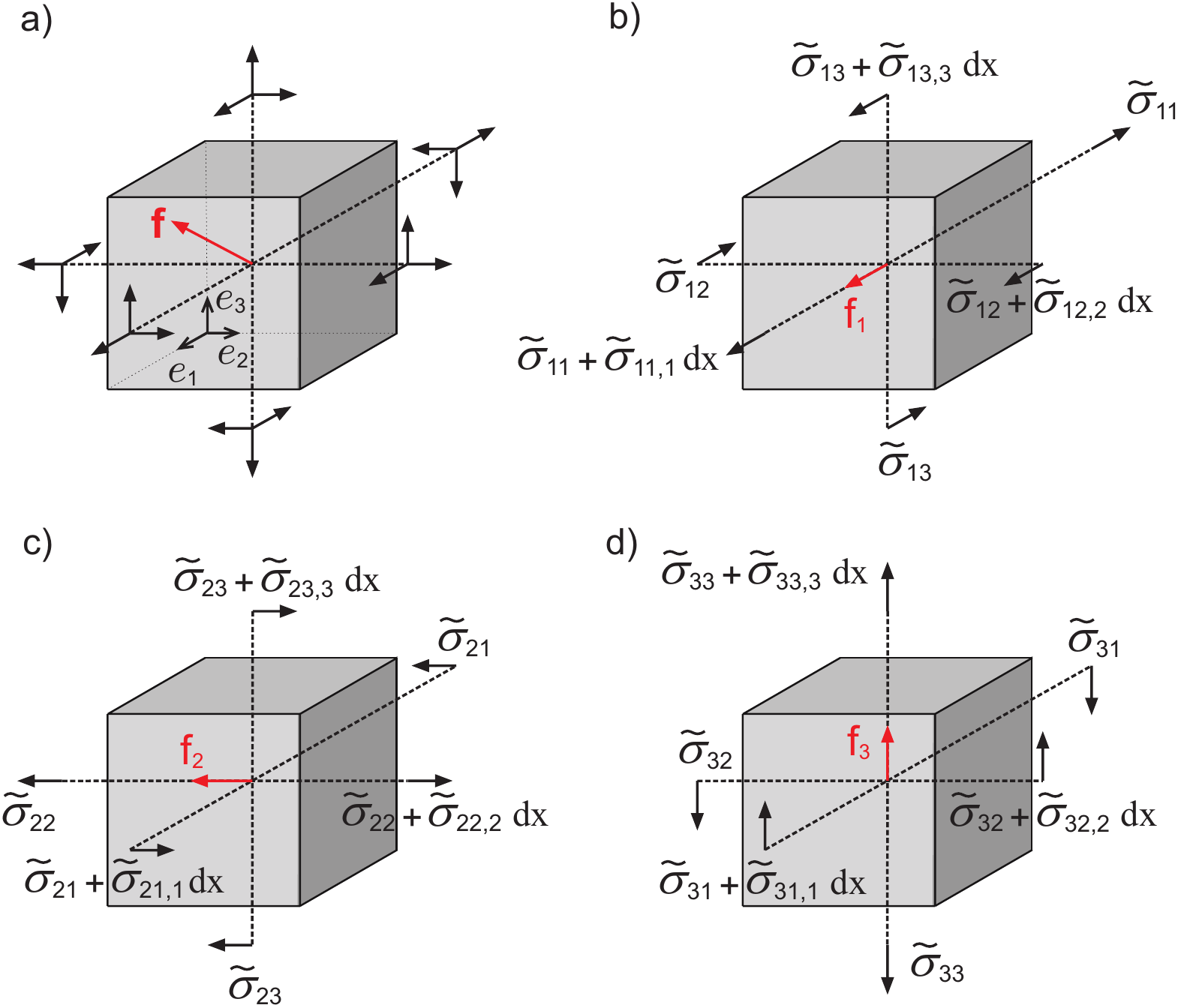}
\caption{a) Cauchy cube $\B_c$ with components of tractions and body force $\f f$. b) Equilibrium of forces in 1-direction. c) Equilibrium of forces in 2-direction. d) Equilibrium of forces in 3-direction.}
\label{PicImpulsBilanzClassic}
\end{figure}

The classical balance of linear momentum considers tractions from the stress tensor $ \widetilde{\fg \sigma}$ as vectors in the center of corresponding faces of the Cauchy cube. Therefore, components of tractions can be simply represented by single arrows in the center of faces as sketched in Fig. \ref{PicImpulsBilanzClassic}a. Gradients of components tangential to the  corresponding faces are disregarded such that the derivation of stress components appear in normal direction of faces only, see Fig. \ref{PicImpulsBilanzClassic}b-d. Further, one assumes the body force $\f f$ to be constant and barycentric within $\B_c$.
The balance of forces sum up the volume integral of the body force $\f f$ and the surface integral of tractions $\widetilde{\fg \sigma}.\f n_i$ of each face $\partial \B_i$ in the direction $\f e_k$. Since the volume is given by $\di V = \di x \, \di x \, \di x$ and each face has the area $\di A = \di x \, \di x$, one obtains three equations from each spatial direction $\f e_k$:
\begin{align}\label{3LinMomEqu1}
&\sum_{i=1}^{6} \int_{\partial B_i} (\widetilde{\fg \sigma}. \f n_i) \cdot \f e_1 \, \di A + \int_{\B} \f f \cdot \f e_1 \, \di V = 0 \notag \\
&\Leftrightarrow (\sigmas_{11} + \sigmas_{11,1} \di x - \sigmas_{11} + \sigmas_{12} + \sigmas_{12,2} \di x - \sigmas_{12} + \sigmas_{13} + \sigmas_{13,3} \di x - \sigmas_{13}) \, \di x \, \di x + {\rm f}_1 \, \di x \, \di x \, \di x = 0 \notag \\
& \Leftrightarrow \sigmas_{11,1} +\sigmas_{12,2} +\sigmas_{13,3} + {\rm f}_1 = 0 \,,
\end{align}
\begin{align}\label{3LinMomEqu2}
&\sum_{i=1}^{6} \int_{\partial B_i} (\widetilde{\fg \sigma} \cdot \f n_i). \f e_2 \, \di A + \int_{\B} \f f \cdot \f e_2 \, \di V = 0 \notag \\
&\Leftrightarrow (\sigmas_{21} + \sigmas_{21,1} \di x - \sigmas_{21} + \sigmas_{22} + \sigmas_{22,2} \di x - \sigmas_{22} + \sigmas_{23} + \sigmas_{23,3} \di x - \sigmas_{23}) \, \di x \, \di x + {\rm f}_2 \, \di x \, \di x \, \di x = 0 \notag \\
& \Leftrightarrow \sigmas_{21,1} +\sigmas_{22,2} +\sigmas_{23,3} + {\rm f}_2 = 0 \,,
\end{align}
\begin{align}\label{3LinMomEqu3}
&\sum_{i=1}^{6} \int_{\partial B_i} (\widetilde{\fg \sigma} \cdot \f n_i). \f e_3 \, \di A + \int_{\B} \f f \cdot \f e_3 \, \di V = 0 \notag \\
&\Leftrightarrow (\sigmas_{31} + \sigmas_{31,1} \di x - \sigmas_{31} + \sigmas_{32} + \sigmas_{32,2} \di x - \sigmas_{32} + \sigmas_{33} + \sigmas_{33,3} \di x - \sigmas_{33}) \, \di x \, \di x + {\rm f}_3 \, \di x \, \di x \, \di x = 0 \notag \\
& \Leftrightarrow \sigmas_{31,1} +\sigmas_{32,2} +\sigmas_{33,3} + {\rm f}_3 = 0 \,,
\end{align}
reading
\begin{align}\label{3LinMomEqu4}
\left(
  \begin{array}{c}
    \sigmas_{11,1} +\sigmas_{12,2} +\sigmas_{13,3} \\
    \sigmas_{21,1} +\sigmas_{22,2} +\sigmas_{23,3} \\
    \sigmas_{31,1} +\sigmas_{32,2} +\sigmas_{33,3} \\
  \end{array}
\right)
+
\left(
  \begin{array}{c}
    {\rm f}_1 \\
    {\rm f}_2 \\
    {\rm f}_3 \\
  \end{array}
\right)
= 0
\quad \Leftrightarrow \quad
\Div \widetilde{\fg \sigma} + \f f = 0 \,,
\end{align}
as a vector equation with help of the divergence operator.
\subsection{Position independency of force couples in rigid bodies}\label{KapPosIndependentM}
\begin{figure}
\centering
\includegraphics[height=80mm]{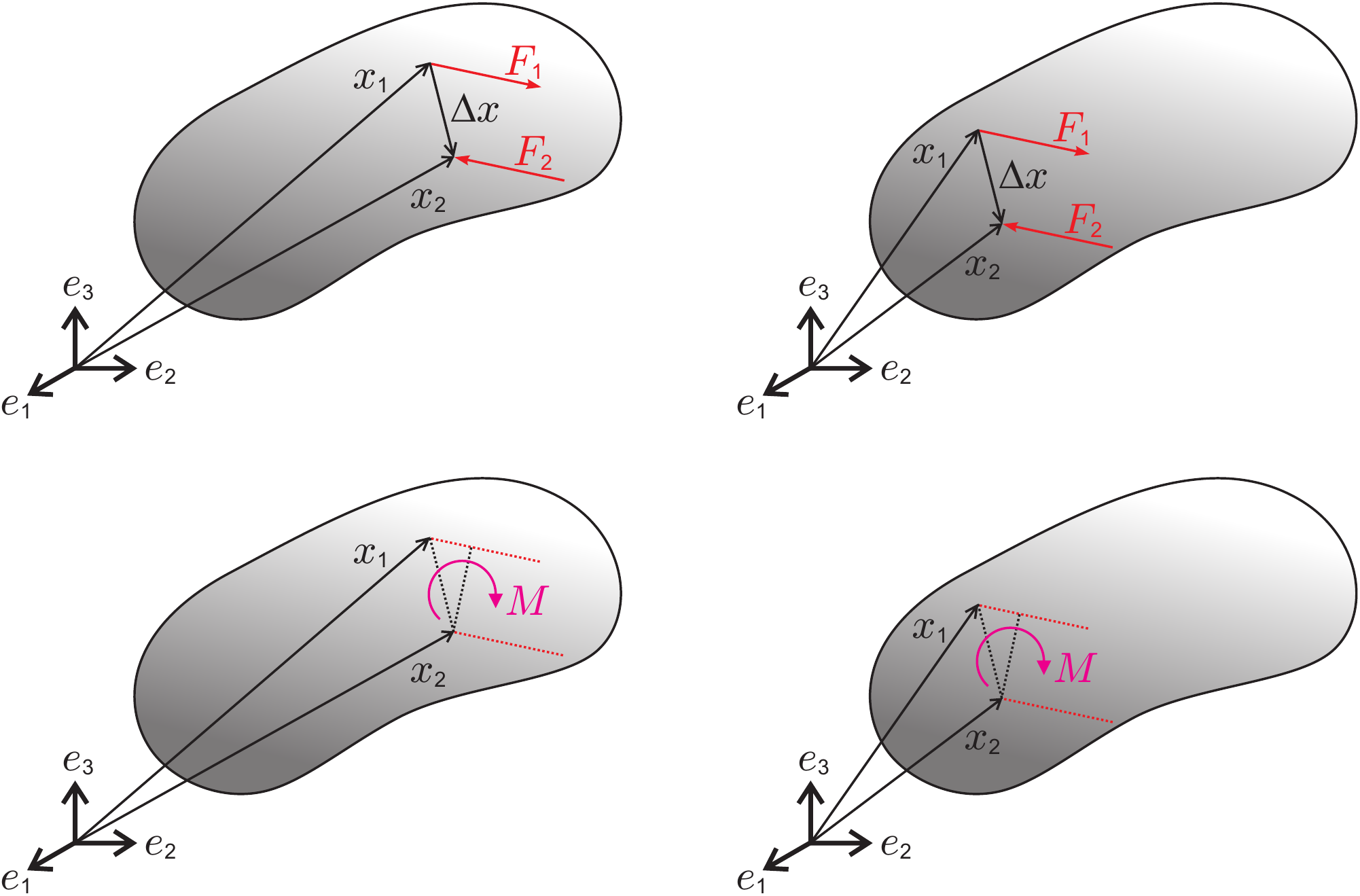}
\caption{A couple of forces applied to a rigid body can be represented by the moment $M$. The spin given to the rigid body does not depend on the position of the moment.}
\label{FreeCouple}
\end{figure}

Let us consider two forces $\f F_1$ and $\f F_2$ applying at position $\f x_1$ and $\f
x_2$, respectively. We presume the following properties:
\begin{align}\label{SumTwoForces}
\f F_1=- \f F_2 \quad , \quad \f F_1,\f F_2 \in \Euklid
\end{align}
and
\begin{align}\label{deltaX}
 \Delta \f x = \f x_2 - \f x_1 \quad , \quad \Delta \f x \neq 0 \quad , \quad \f x_1,\f x_2,\Delta \f x \in \Euklid .
\end{align}
The couple of forces define the moment $\f M$ by their distance $\Delta \f x$ within a cross product
\begin{align}\label{MomTwoForces}
 \f M \colonequals \f x_1 \times \f F_1 + \f x_2 \times \f F_2 =  \f x_1 \times (-\f F_2) + (\f x_1 + \Delta \f x) \times (\f F_2)
 = \Delta \f x \times \f F_2 .
\end{align}
Applying $\f F_1$ and $\f F_2$ onto a rigid body creates no acceleration to its center of mass due to eq.\eqref{SumTwoForces}. However, the forces give spin to the body due to the amount of $\f M$ in eq.\eqref{MomTwoForces}. The spin is generally independent of the
position of the moment $\f M$ within the {\bf rigid body}. This can be shown by
translating the couple of forces by an arbitrary distance $\f x$. The distributive
property of the cross product implies
\begin{align}\label{CoupleTrans}
\f M = (\f x_1+\f x) \times \f F_1 + (\f x_2+\f x) \times \f F_2 = \f x_1 \times \f F_1 + \f x_2 \times \f F_2 + \f x \times \underbrace{(\f F_1+\f F_2)}_{=0} = \Delta \f x \times \f F_2 \,.
\end{align}
Thus, a pure moment is a free vector in space. Vice versa, applying the forces $\f F_1$, $\f F_2$ as couple, they are free vectors as well, see Fig.\eqref{FreeCouple}.
Of course, the cross product is not associative. The cross product of $\Delta \f x$ with
$M$ from eq.\eqref{MomTwoForces} does not vanish
\begin{align}\label{MNotAssoci}
\Delta \f x \times \f M = \Delta \f x \times (\Delta \f x \times \f F_2) \neq
\underbrace{(\Delta \f x \times \Delta \f x)}_{=0} \times \f F_2 \,.
\end{align}
Similarly, the cross product of a position vector $\f x$ with a linear independent force $\f F$ results in a moment of force, which is not linear dependent on $\f x$, thus
\begin{align}\label{xxF}
(\f x \times \f F) \neq 0 \quad \Leftrightarrow \quad \f x \times (\f x \times \f F) \neq 0 \,.
\end{align}
We suppose, that disregarding eq.\eqref{MNotAssoci} or eq.\eqref{xxF} lead Yang et al.~\cite{Yang02} to eq.\eqref{ThirdBalance}, which represents an artificial balance
law\footnote{\cite[p.~25]{hadjesfandiari2014evo} write with regard to the development in
Yang et al.~\cite{Yang02}: ``The symmetric character of the couple stress tensor
[$\widetilde{m}$] is based on an artificial fundamental law for equilibrium of couples,
which has no physical reality."}.
\subsection{Divergence theorem including a cross product on tensorial quantities}\label{DivTheoCrossProd}
Let us consider a second order tensor $\f A=\f a \otimes \f b$, defined by vectors $\f a , \f b \in \Euklid$. On a surface ${\partial V}$ the tensor $\f A$ transforms a normal vector to the surface, $\f n$, into a vector $\f A \cdot\f n$ given by
\begin{align}\label{TransAn}
 \f A \cdot\f n = (\f a \otimes \f b)\, \f n
 = a_i \, b_j \, n_k \, \f e_i \, \Scal{\f e_j, \f e_k} = a_i \, b_j \, n_k \, \delta_{jk} \, \f e_i  = a_i \, b_j \, n_j \, \f e_i = \f a \Scal{\f b,\f n} .
\end{align}
For index notation we use orthogonal unit vectors $\f e_i$ and Einstein's summation convention for repeating subscripts. With help of the divergence theorem the surface integral
\begin{align}\label{IntxcrossA}
\int_{\partial V} \f x \times \f A \cdot \f n  \, {\rm d}{a} =  \int_{\partial V} \f x \times (\f a \otimes \f b) \cdot \f n \, {\rm d}{a}  =  \int_{\partial V} (\f x \times \f a \otimes \f b) \cdot \f n  \, {\rm d}{a} = \int_{V} \Div{(\f x \times \f a \otimes \f b)}  \, {\rm d}{v}
\end{align}
becomes a body integral since vector products are associative. Within eq.\eqref{IntxcrossA} one can express the cross product by using the Levi-Civita tensor
$\fg \epsilon$ by
\begin{align}\label{xtimesA}
\f x \times \f a \otimes \f b =\epsilon_{mni} \, x_m \, a_n \, \f e_i \otimes b_j \, \f e_j = \fg \epsilon : \f x \otimes \f A.
\end{align}
Applying the divergence operator to the expression in eq.\eqref{xtimesA} yields
\begin{align}\label{DivxcrossAResult}
\Div{(\fg  \epsilon : \f x \otimes \f A)}=\frac{\partial \epsilon_{mni} \, x_m \, A_{nj}}{\partial x_k} \, \f e_i \, \Scal{\f e_j,\f e_k} = (\epsilon_{mni} \, \delta_{mk} \, A_{nj} + \epsilon_{mni} \, \, x_m \, A_{nj,k}) \, \delta_{jk} \, \f e_i\\ \notag =(\epsilon_{jni} \, A_{nj} + \epsilon_{mni} \, \, x_m \, A_{nj,j}) \, \f e_i = - \fg \epsilon : \f A + \f x \times \Div{\f A} = 2 \, \axl{\f A} + \f x \times \Div{\f A}.
\end{align}
Thus, the divergence theorem including a cross product on tensorial quantities reads
\begin{align}\label{DivTheoxcrossA}
\int_{\partial V} \f x \times \f A \cdot \f n  \, \di a = \int_{V} 2 \, \axl{\f A} + \f x \times \Div{\f A} \, \di V.
\end{align}
We make use of eq.\eqref{DivTheoxcrossA} for transformations within eq.\eqref{CauchyAndDiv3} and eq.\eqref{CauchyAndDiv4}, respectively.
\end{document}